\documentclass[12pt]{article}
\usepackage{amsmath,theorem,
amssymb,amscd,amsfonts}
\usepackage[all]{xy}
\input xypic

 \newtheorem{corollary}[equation]{Corollary} \newtheorem{lemma}[equation]{Lemma} \newtheorem{proposition}[equation]{Proposition} \newtheorem{theorem}[equation]{Theorem} 
 
\theoremstyle{plain}
\theorembodyfont{\rmfamily}  \newtheorem{definition}[equation]{Definition}  \newtheorem{construction}[equation]{Construction} \newtheorem{remark}[equation]{Remark}\newenvironment{pf}{\noindent \textbf{Proof.}}{\hfill{$\square$}\\} \numberwithin{equation}{section}

\DeclareMathOperator{\Pic}{Pic}
\DeclareMathOperator{\trdeg}{tr.deg}
\DeclareMathOperator{\sh}{sh}
\DeclareMathOperator{\cHom}{\mathcal{H}\mathnormal{om}}
\DeclareMathOperator{\Supp}{Supp}

\newcommand{\oX}{\overline{X}}
\newcommand{\oZ}{\overline{Z}}

\def\st{\,|\,}

\def\discrep{\operatorname{discrep}}

\def\w{\omega}
\def\O{\mathcal{O}}

\def\P{\mathbb{P}}
\def\N{\mathbb{N}}

\def\Z{\mathbb{Z}}
\def\Q{\mathbb{Q}}

\def\a{\alpha}
\def\b{\beta}

\def\D{\Delta}

\def\l{\lambda}
\def\O{\mathcal{O}}

\def\t{\tau}
\def\s{\sigma}

\def\z{\zeta}
\def\r{\rho}

\def\la{\langle} 
\def\ra{\rangle}
\def\rtar{\rightarrow}

\def\det{\operatorname{det}}

\def\rad{\operatorname{rad}}

\def\Stab{\operatorname{Stab}}

\def\Mod{\operatorname{Mod}}
\def\refl{\operatorname{ref}}

\def\ker{\operatorname{ker}}
\def\Hom{\operatorname{Hom}}

\def\End{\operatorname{End}}
\def\Pic{\operatorname{Pic}}
\def\Spec{\operatorname{Spec}}
\def\Proj{\operatorname{Proj}}

\def\Br{\operatorname{Br}}

\def\GL{\operatorname{GL}}
\def\PGL{\operatorname{PGL}}

\def\dim{\operatorname{dim}}

\begin{document}
\title{Canonical singularities of orders over surfaces}
\author{Daniel Chan, Paul Hacking, and Colin Ingalls}
\maketitle

\begin{abstract}
We define and study canonical singularities of orders over surfaces. 
These are noncommutative analogues of Kleinian singularities which arise naturally in the minimal model program for orders over surfaces \cite{CI}. 
We classify canonical singularities of orders using their minimal resolutions (which we define).  
We describe them explicitly as invariant rings for the action of a finite group on a full matrix algebra over a regular local ring. 
We also prove that they are Gorenstein, describe their Auslander-Reiten quivers, and note a simple version of the McKay correspondence.
\end{abstract}

\begin{section}{Introduction}

We define and study canonical singularities of orders over surfaces.
These are noncommutative analogues of Kleinian or Du Val singularities which arise naturally in the minimal model program for orders over surfaces
\cite{CI}.

Let $Z$ be a normal surface over an algebraically closed field $k$ of characteristic $0$. 
An \emph{order} over $Z$ is a coherent sheaf $\O_X$ of $\O_Z$-algebras with generic fibre $k(X):=\O_Z \otimes k(Z)$ a central simple $k(Z)$-algebra.
We also impose a normality condition on $\O_X$ which, if $Z$ is projective, requires that $\O_X$ is \emph{maximal}, that is,
maximal with respect to inclusion among orders in $k(X)$.
We think of the order $\O_X$ as a noncommutative surface $X$ with a finite map to the (commutative) surface $Z=Z(X)$.

Orders provide a way to study division algebras over function fields using the techniques of algebraic geometry.
From this point of view it is natural to study the \emph{birational geometry} of orders.
The modern approach to the birational geometry of varieties is given by Mori's \emph{minimal model program} (see \cite[Sec.~1.4]{KM} for an introduction).
This program was generalised to the case of orders over surfaces in \cite{CI}.

The most important classes of singularities in the minimal model program are the following:
\begin{center}
\{terminal sing.\} $\subset$ \{canonical sing.\} $\subset$ \{log terminal sing.\}
\end{center}
In the case of (commutative) surfaces, the terminal singularities are smooth, the canonical singularities are the Kleinian or Du~Val singularities, and the log terminal singularities are the quotient singularities, i.e., quotients of a smooth surface germ by a finite group. 
The three classes of singularities above are defined using the notion of \emph{discrepancy}. 
In \cite{CI}, terminal, canonical, and log terminal orders were defined by generalising this notion to the case of orders over surfaces 
(see Def.~\ref{definitiondiscrepancy}), and terminal singularities of orders were classified (see Thms.~\ref{tterminal} and \ref{tterminal2}).
In this paper we study canonical and log terminal singularities of orders.
Our main goal is an explicit classification of canonical singularities of orders.

A \emph{birational morphism} $f \colon Y \rightarrow X$ of orders is a birational morphism $Z(Y) \rightarrow Z(X)$ 
of the central surfaces and a morphism $f^{\sharp} \colon f_Z^*\O_X \rightarrow \O_Y$ of $\O_{Z(Y)}$-algebras
which restricts to an isomorphism over the complement of the exceptional locus.
Every order admits a \emph{resolution}, that is, a birational morphism $f \colon Y \rightarrow X$ where $f_Z$ is proper
and the order $Y$ is a terminal order \cite[Cor.~4.5]{CI}.
In particular, the centre $Z(Y)$ is smooth and $\O_Y$ has global dimension $2$.

If $\O_X$ is an order over a surface $Z$ and $P \in Z$ is a general point,
then the fibre $\O_X \otimes k(P)$ is a full matrix algebra over $k$.
This fails in codimension $1$ however. For $C \subset Z$ a curve the structure of the restriction 
$\O_X \otimes k(C)$ determines a cyclic cover of $C$ of some degree $e_C$ (which is nontrivial for finitely many $C$), see Thm.~\ref{cycliccover}.
We refer to this data as the \emph{ramification data} of the order $X$.
We define the \emph{canonical divisor} of $X$ by $K_X=K_Z+\Delta$, where $\Delta = \sum_{C \subset Z} (1-\frac{1}{e_C})C$. 
This is motivated by a formula for the dualising sheaf $\omega_X =\cHom(\O_X,\omega_Z)$ of $X$, see Prop.~\ref{candiv}.

We define the minimal model program for terminal orders over surfaces using the minimal model program for
the surface-divisor pair $(Z,\Delta)$ (see e.g. \cite[Thm.~3.47]{KM}). 
That is, we inductively contract curves $E \subset Z$ such $(K_Z+\Delta) \cdot E < 0$.
At each stage, given a contraction $Z \rightarrow Z'$, we define the order $\O_{X'}$ over $Z'$ as the reflexive hull of the pushforward of $\O_X$.
In particular, applying the minimal model program in the relative context to a resolution $f \colon Y \rightarrow X$,
we obtain a minimal resolution $f' \colon Y' \rightarrow X$ which is uniquely determined up to Morita equivalence, see Cor.~\ref{minresorder}. 
Here we say $f'$ is \emph{minimal} if $K_{Y'} \cdot E  \ge 0$ for every exceptional curve $E$ of $f'_Z$.
We can therefore study singularities of orders in terms of their minimal resolutions. 
This is a natural and powerful approach, as in the commutative setting.

We classify canonical singularities of orders as follows.
We first determine their minimal resolutions, and deduce the classification of the ramification data, see Thm.~\ref{canramdata}.
We then give explicit descriptions of canonical singularities of orders as \emph{quotient singularities}, i.e., invariant rings $B^G$
where $B$ is a full matrix algebra over a regular local ring and $G$ is a finite group acting on $B$, see Sec.~\ref{exp}.
We also show that the log terminal singularities of orders are precisely the quotient singularities (Thm.~\ref{tltisfrt}),
generalising the commutative result mentioned above.

There are three interesting byproducts of our explicit quotient constructions of canonical singularities of orders. 
First, we show that every canonical order is Gorenstein, that is, the dualising sheaf $\omega_X$ is locally isomorphic to $\O_X$ 
as a left and right $\O_X$-module. 
(Note that the Gorenstein property cannot be deduced from the ramification data, so is quite subtle, see Sec.~\ref{ltgor}.)
Second, in studying the normality condition, we are led to consider certain reflexive 
modules called permissible modules. We observe that the number of permissible modules equals 
the number of exceptional curves in the minimal resolution plus one.
This suggests some version of the McKay correspondence holds for canonical orders. 
Finally, we compute the Auslander--Reiten quivers of canonical singularities of orders. 
This relates the classifications of Artin \cite{A86} and Reiten--Van den Bergh \cite{RVdB} in the case of canonical orders. 

Canonical orders can be applied to the study of $3$-fold conic bundles \cite{CHI} and the classification of maximal orders over projective surfaces 
\cite{CK}.

\medskip
\noindent
\textbf{Notation}:
Throughout this paper, we work over an algebraically closed field $k$ of characteristic $0$.
By a \emph{variety}, we mean a Noetherian integral separated scheme over $k$. 
A \emph{curve} (resp. \emph{surface}) is a variety of dimension $1$ (resp. $2$).

A \emph{$\Q$-divisor} $D$ on a normal variety $Z$ is a finite formal linear combination of codimension $1$ subvarieties
with rational coefficients. 
A $\Q$-divisor $D$ is \emph{$\Q$-Cartier} if $nD$ is Cartier for some $n \in \N$.
A $\Q$-Cartier $\Q$-divisor $D$ on a variety $Z$ is \emph{nef} (``numerically eventually free")
if $D \cdot C \ge 0$ for every curve $C \subset Z$. Similarly, if $g \colon Z \rightarrow W$ a morphism, 
we say $D$ is \emph{$g$-nef}, or \emph{nef over W}, if $D \cdot C \ge 0$ for every curve $C$ contracted by $g$.

\medskip
\noindent
\textbf{Acknowledgements}: We thank Bill Crawley-Boevey
and Alastair King for helpful information about Auslander-Reiten
quivers. The first author was supported by an Australian Research Council grant.
The second author was partially supported by NSF grant DMS-0650052.

\end{section}

\begin{section}{Birational geometry of orders}  \label{sprelim}

A version of Mori's minimal model program for orders over surfaces was developed in \cite{CI}. 
In this section, we review the basic notions of this theory with an emphasis on the local aspects. 
In particular, we recall the notion of discrepancy, and use it to define terminal, canonical, and log terminal orders.
We also define minimal resolutions of orders, and prove that they are uniquely determined up to Morita equivalence.

\subsection{Normal orders}

Let $Z$ be a surface with  function field $k(Z)$.  We define an  
{\it order} $\O_X$ over $Z$ to be  a torsion-free coherent sheaf of $\O_Z$-algebras such that
$k(X) := \O_X \otimes k(Z)$ is a central simple $k(Z)$-algebra.
We often refer to the order simply  as $X$ and call $Z$ the {\it centre} of $X$. 
The \emph{degree} of $X$ is the integer $\deg X = \sqrt{\dim_{k(Z)} k(X)}$. 
These orders are our noncommutative analogues of surfaces. 
 
In the theory of orders there are various noncommutative analogues of the class of normal varieties.
The simplest of these is the class of maximal orders, defined as follows: an order $X$ is {\it maximal} if it is 
maximal with respect to inclusion among orders inside $k(X)$. 
When $Z$ is projective, we are primarily interested in maximal orders. 
However, as Artin noted in \cite[Sec.~2]{A81}, maximality is not preserved under \'{e}tale localization on $Z$. 
We define the class of normal orders which includes the maximal orders and is stable under henselisation or completion at a point $P \in Z$.
   
\begin{definition} \label{dnormal}
Let $X$ be an order over a surface $Z$. 
The {\it dualising sheaf} of $X$ is the $\O_X$-bimodule $\w_X := \cHom_Z(\O_X,\w_Z)$.
We say $X$ is \emph{normal} if $Z$ is normal and $\O_X$ satisfies the following conditions. 
\begin{itemize}
\item[($R_1$)] 
For all curves $C \subset Z$
\begin{enumerate}
\item the local ring $\O_{X,C}$ is hereditary and $\w_{X,C}$ is isomorphic to $\O_{X,C}$ as a left and right $\O_{X,C}$-module 
(but not necessarily as a bimodule), and
\item $\O_{X,C}$ is maximal whenever $k(C)$ has finite transcendence degree over $k$.  
\end{enumerate}
\item[($S_2$)] $\O_X$ is a reflexive $\O_Z$-module.
\end{itemize}
\end{definition}
Condition $(R_1)(1)$ is equivalent to the following: after an 
\'{e}tale base change $\O_{Z,C} \rightarrow R$, the algebra $\O_{X,C} \otimes R$ has the standard form
$$\left(\begin{array}{cccc} 
                                R    & \cdots & \cdots  & R\\
                                tR   & \ddots &         & \vdots \\ 
				\vdots & \ddots & \ddots  & \vdots \\
				tR   & \cdots & tR    & R
                               \end{array} \right)^{f \times f},$$
where $t \in R$ is a local parameter, the displayed matrix is $e \times e$, and  $ef=\deg X$. 
To clarify condition $(R_1)(2)$, note that the surfaces $Z$ we consider are often the source of a proper birational morphism with target the spectrum of a complete local ring of dimension $2$. 
In this case $(R_1)(2)$ requires maximality at the generic points of the exceptional curves. 
We only consider normal orders in what follows.

\subsection{Ramification and the canonical divisor} \label{subsectionramification}

The following well-known result leads to the concept of ramification for orders over surfaces. 
This is the key tool which allows us to study orders geometrically. 

\begin{theorem}\cite[Prop.~2.4]{CI} \label{cycliccover}
Let $\O_X$ be a normal order over a surface $Z$.  
Let $C \subset Z$ be a curve and $J$ the radical of $\O_{X,C}$.
Then $\O_{X,C}/J \simeq \prod_{i=1}^r L^{n \times n}$, 
where $L$ is a cyclic field extension of $k(C)$.
Moreover, $\O_{X,C}$ is maximal iff $r=1$.
Finally, $L$ is determined by the Brauer class of $k(X)$.  
\end{theorem}
For $C \subset Z$ a curve, the {\it ramification index} $e_C$ of $X$ over $C$ is the degree of the extension $L^r/k(C)$.
The {\it discriminant} is the union of the ramification curves 
$$D := \bigcup_{e_C > 1} C.$$
The extension $L^r/k(C)$ determines a smooth cover $\tilde{C}$ of the normalisation of $C$ of degree $e_C$
(where $\tilde{C}$ is a disjoint union of $r$ isomorphic curves).  
We let $\tilde{D}$ denote the disjoint union of the covers $\tilde{C}$ as $C$ varies over the ramification 
curves. We define the {\it ramification data} of $X$ to be the centre $Z$ of $X$, the discriminant $D \subset Z$, 
and the cover $\tilde{D} \rightarrow D$. We define the ramification data of a Brauer class $\alpha \in \Br(k(Z))$ 
over a curve $C \subset Z$ to be that of any maximal order $\O_X$ such that the Brauer class of $k(X)$ equals $\alpha$.

\begin{remark} \label{blowupram}
Let $X$ be a normal order over a surface $Z$ and $\alpha$ the Brauer class of $k(X)$.
Let $\pi \colon \tilde{Z} \rightarrow Z$ be the blowup of a smooth point $P \in Z$ and $E \simeq \P^1$ the exceptional curve.
Then the ramification of $\alpha$ over $E$ can be explicitly determined from the ramification data of $X$ at $P \in Z$.
See \cite[Lem.~3.4]{CI}.
\end{remark}
 

We define the canonical divisor of a normal order $X$ using its ramification data. With notation as above, define the \emph{ramification divisor}
$$\Delta = \sum_{C} \left( 1 - \frac{1}{e_C} \right) C,$$
a $\Q$-divisor on $Z$.
We define the {\it canonical divisor} of $X$ by $K_X = K_Z + \D$.
This is motivated by the following.

\begin{proposition} \label{candiv}\cite[Prop.~5(2)]{CK}
Let $X$ be a normal order of degree $m$ over a surface $Z$. Then, in codimension $1$, there is
a natural isomorphism of $\O_X$-bimodules
$$\w_X^{\otimes m} \simeq \O_X \otimes_{\O_Z} \O_Z(m(K_Z + \D)).$$
\end{proposition}

\subsection{The minimal model program relative to a Brauer class} \label{MMPtriples}

Let $Z$ be a normal surface, $\alpha \in \Br k(Z)$ a Brauer class, and $\Delta$ the ramification divisor
of a normal order $X$ over $Z$ such that the Brauer class of $k(X)$ equals $\alpha$.
(Note that if $\O_X$ is maximal then $\Delta$ is determined by $\alpha$, but we do not want to assume this.)
We refer to the data $(Z,\Delta,\alpha)$ as a \emph{triple}.
If the order $X$ is given, we call $(Z,\Delta,\alpha)$ the \emph{associated triple}.
In this section we describe the birational geometry of triples.
In Sec.~\ref{MMPorders} we deduce the analogous results for orders.

\begin{remark}
The ramification data of $X$ is determined by the associated triple.
Indeed, for $C \subset Z$ a curve, the ramification index $e_C$ is determined by the coefficient of $C$ in $\Delta$
and a connected component of the cyclic cover $\tilde{C} \rightarrow C$ is determined by $\alpha$.
\end{remark}
  
If $(Z,\Delta,\alpha)$ is a triple and $f \colon W \rightarrow Z$ is a birational morphism from a normal surface $W$, 
we write $\Delta_W$ (or just $\Delta$ if no confusion seems likely) for the $\Q$-divisor on $W$ determined by $\Delta$ and $\alpha$ as follows: 
$$
\Delta_W=\Delta'+ \sum \left(1-\frac{1}{e_i}\right)E_i
$$ 
where $\Delta'$ is the strict transform of $\Delta$, the $E_i$ are the exceptional curves, and $e_i$ is the ramification index of
$\alpha$ over $E_i$. 
Similarly, if $g \colon Z \rightarrow V$ is a birational morphism to a normal surface $V$, we write $\Delta_V$ (or just $\Delta$)
for the $\Q$-divisor $g_*\Delta$ on $V$.
If $X$ is a normal order with associated triple $(Z,\Delta,\alpha)$, then $(W,\Delta_W,\alpha)$ is the associated triple of a normal order
$\O_Y$ obtained from $f^*\O_X$ by ``saturating" at the exceptional curves, see \ref{blowupconstruction} for details.
Similarly, $(V,\Delta_V,\alpha)$ is the associated triple of the reflexive hull of $g_*\O_X$.

\medskip
A \emph{resolution} of a surface $Z$ is a proper birational morphism $f \colon W \rightarrow Z$ from a smooth surface $W$.
\begin{definition} \label{definitiondiscrepancy}
Let $(Z,\Delta,\alpha)$ be a triple.
Let $f \colon W \rightarrow Z$ be a resolution of $Z$.
Then we have an equality 
$$
K_W+\Delta=f^*(K_Z+\Delta)+ \sum a_iE_i
$$
in $\Pic(W) \otimes_{\Z} \Q$, where the $E_i$ are the exceptional curves and $a_i \in \Q$.
The \emph{discrepancy} $a(E_i,Z,\Delta,\alpha)$ of $E_i$ over $(Z,\Delta,\alpha)$ is the rational number $e_ia_i$,
where $e_i$ is the ramification index of $\alpha$ over $E_i$.

The \emph{discrepancy} $\discrep(Z,\Delta,\alpha)$ of $(Z,\Delta,\alpha)$ is the infimum of the discrepancies
$a(E,Z,\Delta,\alpha)$, where $E$ runs through all exceptional curves of all resolutions of $Z$.
(If the discrepancies are unbounded we write $\discrep(Z,\Delta,\alpha)=-\infty$.)

We say $(Z,\Delta,\alpha)$ is \emph{terminal} if $\discrep(Z,\Delta,\alpha) > 0$, \emph{canonical} if $\discrep(Z,\Delta,\alpha) \ge 0$, 
and \emph{log terminal} if $\discrep(Z,\Delta,\alpha) > -1$.
\end{definition}

We have an \'etale local structure theorem for terminal triples $(Z,\Delta,\alpha)$.

\begin{theorem}\cite{CI} \label{tterminal}
Let $(Z,\Delta,\alpha)$ be a triple.
Then $(Z,\Delta,\alpha)$ is terminal iff $Z$ is smooth and, \'etale locally at each point $P \in Z$, 
the ramification data is of one of the following types.
\begin{enumerate}
\item $D=\emptyset$.
\item $(P \in D)$ is smooth, and $\tilde{D} \rightarrow D$ is \'etale of degree $e$.
\item $(P \in D)$ is a node, the cyclic cover $\tilde{D}$ of
the normalisation of $D$ has degrees $e$ and $ne$ over the two branches of $D$,
and $\tilde{D}$ has ramification index $e$ over $P$ on each branch of $D$.   
\end{enumerate}
\end{theorem}

\begin{remark} \label{blowupterminal}
If $(Z,\Delta,\alpha)$ is a terminal triple and $\pi \colon \tilde{Z} \rightarrow Z$ is the blowup of a point $P \in Z$, then 
$(\tilde{Z},\Delta,\alpha)$ is also terminal. To see this, it is enough to compute the ramification of $\alpha$ over
the exceptional curve $E \simeq \P^1$ in each of the cases of Thm.~\ref{tterminal}.
We find that $\alpha$ is unramified over $E$ in cases (1) and (2), and $\alpha$ has ramification index $e$ over $E$ case (3).
In case (3) the cover $\tilde{E} \rightarrow E$ is the unique degree $e$ cyclic cover which is totally ramified over the $2$ 
intersection points of $E$ with the strict transform of $D$ and \'etale elsewhere. 
\end{remark} 

\begin{definition}
A \emph{resolution} of $(Z,\Delta,\alpha)$ is a resolution $f \colon \tilde{Z} \rightarrow Z$ of $Z$ such that $(\tilde{Z},\Delta,\alpha)$ is terminal.
\end{definition}

Resolutions always exist \cite[Cor.~3.6]{CI}.

We next describe how to compute the discrepancy of a triple from a single resolution.

\begin{lemma}
Let $(Z,\Delta,\alpha)$ be a triple, $f \colon W \rightarrow Z$ a birational morphism, and 
$g \colon V \rightarrow W$ the blowup of a smooth point $P \in W$. 
Let $E_j$ be the exceptional curves of $f$, $a_j=a(E_j,Z,\Delta,\alpha)$, $e_j$ the ramification index of $\alpha$
along $E_j$, and $\nu_j$ the multiplicity of $E_j$ at $P$.
Write $\Delta_W = \sum (1-\frac{1}{e_i})C_i$ and let $\mu_i$ be the multiplicity of $C_i$ at $P$.
Let $E$ be the exceptional curve of $g$ and $e$ the ramification index of $\alpha$ along $E$. Then
$$a(E,W,\Delta,\alpha)=e\left(1 + \left(1-\frac{1}{e}\right) - \sum \mu_i\left(1-\frac{1}{e_i}\right)\right)$$ 
and
$$a(E,Z,\Delta,\alpha)=a(E,W,\Delta,\alpha)+ e\sum \nu_j\frac{a_j}{e_j}.$$
\end{lemma}
\begin{pf}
Let $\Delta_W'$ denote the strict transform of $\Delta_W$ on $V$.
We have $K_V=g^*K_W+E$, $\Delta_V=\Delta_W'+ (1-\frac{1}{e})E$, and $\Delta_W'= g^*\Delta_W - (\sum \mu_i(1-\frac{1}{e_i}))E$.
Hence
$$K_V+\Delta_V=g^*(K_W+\Delta_W)+\left(1+\left(1-\frac{1}{e}\right)-\sum \mu_i\left(1-\frac{1}{e_i}\right)\right)E.$$
Now $K_W+\Delta_W=f^*(K_Z+\Delta)+\sum \frac{a_j}{e_j}E_j$ and $g^*E_j=E_j'+\nu_jE$ where $E_j'$ is the strict transform of $E_j$.
So $$K_V+\Delta_V=g^*(K_W+\Delta_W)+ \frac{1}{e}\cdot a(E,W,\Delta,\alpha)E$$
$$=g^*f^*(K_Z+\Delta)+\sum \frac{a_j}{e_j}E_j'+\left(\frac{1}{e}\cdot a(E,W,\Delta,\alpha)+\sum \nu_j\frac{a_j}{e_j}\right)E.$$
\end{pf}

\begin{lemma} \label{discrepinfty}
Let $(Z,\Delta,\alpha)$ be a triple. Then either $\discrep(Z,\Delta,\alpha) \ge -1$ or $\discrep(Z,\Delta,\alpha)=-\infty$.
\end{lemma}
\begin{pf}
(Cf. \cite[Cor.~2.31(1)]{KM})
Suppose $f \colon W \rightarrow Z$ is a resolution of $Z$ and $E \subset W$ is an exceptional curve such that $a(E,Z,\Delta,\alpha)=-(1+c)$ 
for some $c > 0$.
Let $P_0 \in E$ be a general point, $g_1 \colon W_1 \rightarrow W_0=W$ the blowup of $P_0$, and 
$E_1 \subset W_1$ the exceptional curve.  
Then $\alpha$ is unramified over $E_1$ and $a(E_1,Z,\Delta,\alpha)=-\frac{c}{e}$. 
Given $g_n \colon W_n \rightarrow W_{n-1}$ with exceptional curve $E_n$, 
let $P_n$ be the intersection of the strict transform of $E$ and $E_n$ and 
$g_{n+1} \colon W_{n+1} \rightarrow W_n$ the blowup of $P_n$ with exceptional curve $E_{n+1}$.
Then $\alpha$ is unramified over $E_n$ and $a(E_n,Z,\Delta,\alpha)=-\frac{nc}{e}$ for all $n$.
\end{pf}

\begin{proposition} \label{computediscrep}
Let $(Z,\Delta,\alpha)$ be a triple and $f \colon W \rightarrow Z$ a resolution of $(Z,\Delta,\alpha)$.
Assume that the support of the divisor $\Delta'+\sum E_i$ has normal crossing singularities, where
$\Delta'$ is the strict transform of $\Delta$ and the $E_i$ are the exceptional curves. 
For $C \subset Z$ a curve let $e_C$ denote the ramification index of $(Z,\Delta,\alpha)$ over $C$.
Then $\discrep(Z,\Delta,\alpha) \neq -\infty$ iff $a(E_i,Z,\Delta,\alpha) \ge -1$ for all $i$, and, in this case,
$$\discrep(Z,\Delta,\alpha) = \min \left\{\left\{a(E_i,Z,\Delta,\alpha) \right\}, \left\{\frac{1}{e_C}\right\}, 1\right\}$$
\end{proposition}
\begin{pf}
By definition $\discrep(Z,\Delta,\alpha) \le a(E_i,Z,\Delta,\alpha)$.
If $F$ is the exceptional divisor of the blowup a general point of $C \subset Z$ then $\alpha$ is unramified over $F$ and 
$a(F,Z,\Delta,\alpha)=\frac{1}{e_C}$, so $\discrep(Z,\Delta,\alpha) \le \frac{1}{e_C}$. 
Blowing up a general point of $Z$ gives $\discrep(Z,\Delta,\alpha) \le 1$.

Note that if $(W,\Delta,\alpha)$ is a terminal triple and $W' \rightarrow W$ is a resolution of $W$, then the ramification index of $\alpha$
over an exceptional curve $E$ is equal to the ramification index of $(W,\Delta,\alpha)$ over some curve $C \subset W$.  
Indeed, $W' \rightarrow W$ is a composition of blowups, so this follows from the computation in Rem.~\ref{blowupterminal}.

If $a(E_i,Z,\Delta,\alpha) < -1$ for some $i$ then $\discrep(Z,\Delta,\alpha)=-\infty$ by Lem.~\ref{discrepinfty}.
Suppose now that $a(E_i,Z,\Delta,\alpha) \ge -1$ for all $i$.
Let $E$ be an exceptional divisor of a resolution $W' \rightarrow Z$. Let $e_i$ denote the ramification index of $\alpha$ over $E_i$.
We claim that $a(E,Z,\Delta,\alpha) \ge b$ where
$$b:=\min \left\{ \left\{a(E_i,Z,\Delta,\alpha)\right\}, \left\{\frac{1}{e_i}\right\}, \left\{\frac{1}{e_C}\right\}, 1\right\}$$
We may assume that $W' \rightarrow Z$ factors through $W$. Then the morphism $W' \rightarrow W$ is a 
composition of blowups, and we may assume that $E$ is the exceptional divisor of the last blowup. 
We argue by induction on the number $n$ of blowups.
Let $W' \rightarrow W''$ be the last blowup and $g \colon W'' \rightarrow Z$ the induced morphism.
So $a(F,Z,\Delta,\alpha) \ge b$ for every $g$-exceptional curve $F$ by the induction hypothesis.
Let $P \in W''$ be the point we blowup and $G$ the union of the exceptional curves of $g$ and the support of the strict transform of $\Delta$ 
(which has normal crossing singularities by assumption). We have the following cases.
\begin{enumerate}
\item $P \notin G$. Then $a(E,Z,\Delta,\alpha)=1$.
\item $P$ is a smooth point of $G$. Let $\Gamma$ be the component of $G$ through $P$, $e$ the ramification index of $(W'',\Delta,\alpha)$ over 
$\Gamma$, and $a=a(\Gamma,Z,\Delta,\alpha)$ if $\Gamma$ is exceptional and $a=0$ otherwise. Then $\alpha$ is unramified over $E$ and 
$a(E,Z,\Delta,\alpha)=\frac{1}{e}(1+a)$. Thus $a(E,Z,\Delta,\alpha) \ge \frac{1}{e}$ if $a \ge 0$ and $a(E,Z,\Delta,\alpha) \ge a$ if $a < 0$.
\item $P$ is a node of $G$. Let $\Gamma_1,\Gamma_2$ be the two components of $G$ through $P$, $ne,e$ the ramification indices (we allow $e=1$),
and $a_i=a(\Gamma_i,Z,\Delta,\alpha)$ if $\Gamma_i$ is exceptional and $a_i=0$ otherwise. Then $\alpha$ has ramification index $e$ over $E$
and $a(E,Z,\Delta,\alpha)=e \cdot (\frac{1}{ne}(1+a_1)+\frac{1}{e}a_2) = \frac{1}{n}(1+a_1)+a_2$.
Hence $a(E,Z,\Delta,\alpha) \ge a_2$ because $a_1 \ge b \ge -1$.
Finally, if $a_2=0$, $a(E,Z,\Delta,\alpha) \ge \frac{1}{ne}$ if $a_1 \ge 0$ and
$a(E,Z,\Delta,\alpha) \ge a_1$ if $a_1 < 0$.
\end{enumerate}
We conclude that $a(E,Z,\Delta,\alpha) \ge b$. 

It remains to observe that the terms $\frac{1}{e_i}$ may be omitted in the definition of $b$.
If $b=\frac{1}{e_i}$ for some $i$ then in particular $(Z,\Delta,\alpha)$ is terminal. So $e_i=e_C$ for some $C \subset Z$ by the remark above.
\end{pf}

If $(Z,\Delta,\alpha)$ is a terminal triple and $Z$ is projective, 
we may apply the log minimal model program to the pair $(Z,\Delta)$ \cite[Thm.~3.47]{KM}.
We obtain the following result.

\begin{theorem}\cite[Cor.~3.20]{CI} \label{minmodeltriple}
Let $(Z,\Delta,\alpha)$ be a terminal triple such that $Z$ is projective.
Then there exists a birational morphism $f \colon Z \rightarrow W$ such that $(W,\Delta,\alpha)$ is terminal and one of the following holds.
\begin{enumerate}
\item $K_W+\Delta$ is nef.
\item $W$ is the total space of a $\P^1$-bundle over a curve and $(K_W+\Delta) \cdot F < 0$ where $F$ is a fibre.
\item $W$ is isomorphic to $\P^2$ and $-(K_W+\Delta)$ is ample.
\end{enumerate}
Moreover, in case (1) $W$ is uniquely determined.
\end{theorem}
\begin{pf}
We refer to \cite{CI} for the full proof.
Here we include a conceptual proof that $(W,\Delta,\alpha)$ is terminal, cf. \cite[Cor.~3.43]{KM}.
By construction, the morphism $f$ is a composition of birational morphisms
$$Z=Z_0 \stackrel{f_1}{\rightarrow} Z_1 \stackrel{f_2}{\rightarrow} \cdots \stackrel{f_n}{\rightarrow} Z_n=W$$
where $f_i \colon Z_{i-1} \rightarrow Z_{i}$ contracts a single curve $E_i$, and $(K_{Z_{i-1}}+\Delta)\cdot E_i < 0$.
So it suffices to show that if $(Z,\Delta,\alpha)$ is terminal and $f \colon Z \rightarrow W$ is a birational morphism
with a unique exceptional curve $E$ such that $(K_Z+\Delta) \cdot E < 0$, then $(W,\Delta,\alpha)$ is terminal. Write
$$K_Z+\Delta=f^*(K_W+\Delta)+aE,$$
then $(K_Z+\Delta) \cdot E =aE^2< 0$ and $E^2 <0$, so $a>0$.
Now let $g \colon V \rightarrow Z$ be a resolution.
Then
$$K_V+\Delta=g^*(K_Z+\Delta)+\sum a_jF_j$$
where the $F_j$ are the exceptional curves and $a_j > 0$ for all $j$ because $(Z,\Delta,\alpha)$ is terminal.
Write $g^*E=E'+\sum \mu_jF_j$, where $E'$ is the strict transform of $E$ and $\mu_{j} \ge 0$ for all $j$.
Then
$$K_V+\Delta=(fg)^*(K_W+\Delta)+aE'+\sum_j (a_j+a\mu_j)G_j$$
Hence $(W,\Delta,\alpha)$ is terminal.

Uniqueness in case (1) follows from \cite[Thm.~3.52(2)]{KM}.
\end{pf}

\begin{definition}
A \emph{minimal resolution}  of a triple $(Z,\Delta,\alpha)$ is a resolution $f \colon \tilde{Z} \rightarrow Z$ of $Z$ such that 
$K_{\tilde{Z}}+\Delta$ is $f$-nef, i.e., $(K_{\tilde{Z}}+ \Delta) \cdot E \ge 0$ for every exceptional curve $E$.
\end{definition}

\begin{theorem} \label{minrestriple}
A triple $(Z, \Delta, \alpha)$ has a unique minimal resolution.
\end{theorem}
\begin{pf}
Let $f \colon W \rightarrow Z$ be a resolution of $(Z,\Delta,\alpha)$.
We apply the relative log minimal model program to the pair $(W,\Delta)$ over $Z$.
That is, we apply the log minimal model program to $(W,\Delta)$ as above but we only
contract curves $E \subset W$ which are contracted by $f$, see \cite[Sec.~3.31]{KM}.
We obtain a birational morphism $W \rightarrow V$ over $Z$ such that $(V,\Delta,\alpha)$ is terminal
and $K_V+\Delta$ is nef over $W$, i.e., $V \rightarrow Z$ is a minimal resolution of $(Z,\Delta,\alpha)$.
Uniqueness follows from \cite[Thm.~3.52(2)]{KM}.
\end{pf}

\begin{remark} \label{minresremark}
Let $(Z,\Delta,\alpha)$ be a triple and $f \colon \tilde{Z} \rightarrow Z$ a resolution of $(Z,\Delta,\alpha)$, 
and write
$$K_{\tilde{Z}}+\Delta=f^*(K_Z+\Delta)+\sum a_iE_i$$
where the $E_i$ are the exceptional curves.
If $f$ is minimal then $a_i \le 0$ for every $i$ by \cite[Lem.~3.41]{KM}.
\end{remark}

\subsection{The minimal model program for orders} \label{MMPorders}

We now describe the corresponding results for orders.
Let $X$ be a normal order over a surface with associated triple $(Z,\Delta,\alpha)$.
We say $X$ is \emph{terminal}, \emph{canonical}, or \emph{log terminal} if $(Z,\Delta,\alpha)$ is so.

\begin{theorem} \cite{CI} \label{tterminal2}
Let $X$ be a terminal order over a surface $Z$ and $P \in Z$ a closed point.
Then, in the cases $(1),(2),(3)$ of Thm.~\ref{tterminal}, the order $\O_X \otimes \hat{\O}_{Z,P}$ 
is isomorphic to a full matrix algebra over the ring below.
Let $R = k[[u,v]]$ and $S = R\langle x,y \rangle/(x^e-u,y^e-v,yx-\zeta xy)$ where $\zeta$ is a primitive $e$-th root of unity.
\begin{enumerate}
\item \quad $R$
\item $$\begin{pmatrix}
R & \cdots & \cdots & R \\
uR & \ddots & \ddots & \vdots \\
\vdots & \ddots & \ddots & \vdots \\
uR & \cdots & uR & R
\end{pmatrix} \subset R^{e \times e}$$
\item $$\begin{pmatrix}
S & \cdots & \cdots & S \\
xS & \ddots & \ddots & \vdots \\
\vdots & \ddots & \ddots & \vdots \\
xS & \cdots & xS & S
\end{pmatrix} \subset S^{n \times n}$$
\end{enumerate}
The centre is $R$ in each case, and the order ramifies over $(u=0)$ with index $e$ in (2), 
and over $(u=0)$ with index $ne$ and $(v=0)$ with index $e$ in (3).

In particular, a terminal order has global dimension $2$ locally over its centre.
\end{theorem}

For an order $X$, let $\Mod X$ denote the category of right $\O_X$-modules.
We say orders $X_1$ and $X_2$ are \emph{Morita equivalent}, and write $X_1 \simeq_M X_2$, 
if there is an equivalence of categories $\Mod X_1 \simeq \Mod X_2$.
Note that the equivalence induces an isomorphism $Z_1 \simeq Z_2$ of the centres \cite[Prop.~6.8]{AZ}.
We say $X_1$ and $X_2$ are \emph{Morita equivalent in A} if there are identifications $k(X_1)=k(X_2)=A$
and a compatible equivalence $\Mod X_1 \simeq \Mod X_2$.

\begin{corollary}\label{terminalMoritaequiv}
Let $X_1$ and $X_2$ be terminal orders over a surface $Z$ with $k(X_1)=k(X_2)=A$.
Assume that $\O_{X_1,C}=\O_{X_2,C} \subset A$ for every curve $C \subset Z$ such that $k(C)$ is not
of finite transcendence degree over $k$.
Then $X_1$ and $X_2$ are Morita equivalent in $A$.
\end{corollary}
\begin{pf}
We follow the proof of the projective case \cite[Cor~2.13]{CI}. 
Consider the reflexive hull $P$ of the $(\O_{X_1},\O_{X_2})$-bimodule $\O_{X_1}\O_{X_2} \subset A$. 
Note that since the orders $\O_{X_i}$ are terminal they have global dimension $2$ locally over the centre.
So $P$ is locally projective on both sides by the Auslander-Buchsbaum formula.
We claim that the natural maps
$$ \O_{X_1} \rtar \End P_{X_2}  \hspace{1cm} \O_{X_2} \rtar \End {}_{X_1} P$$
are isomorphisms. All the objects are reflexive so we need only check this locally at codimension one primes. 
It holds at curves $C$ such that $\trdeg_k k(C) < \infty$ since $\O_{X_1,C},\O_{X_2,C}$ are maximal 
and at the remaining curves by our assumption.
So $P$ gives a Morita equivalence $X_1 \simeq_M X_2$.
\end{pf}

\begin{definition}\label{birationalmorphism} 
A {\it birational morphism} $f: Y \rightarrow X$ of normal orders over surfaces 
is a pair $f=(f_Z,f^\#)$ where $f_Z : Z(Y) \rightarrow Z(X)$ is a birational morphism and 
$f^\# : f_Z^*\O_X \rightarrow \O_Y$ is an $\O_{Z(Y)}$-algebra map which restricts to an isomorphism
over the complement of the exceptional curves. We say $f$ is \emph{proper} if $f_Z$ is so.
\end{definition}

\begin{construction} \label{blowupconstruction}
Let $X$ be a normal order over a surface $Z$.
Let $g \colon W \rightarrow Z$ be a birational morphism.
For each exceptional curve $E$ of $g$, let $\Lambda_E$ be a maximal order in $k(X)$ over $\O_{W,E}$
containing $\O_{X,g(E)}$. 
Let $\O_Y$ be the unique normal order in $k(X)$ over $W$ such that $\O_{Y,E}=\Lambda_E$ for $E$ an exceptional
curve and $\O_{Y,C}=\O_{X,g(C)}$ for $C$ a non-exceptional curve. 
Then the natural $\O_W$-algebra map $g^*\O_X \rightarrow \O_Y$
defines a birational morphism $f \colon Y \rightarrow X$ with $f_Z=g$.

Note that $Y$ is \emph{not} uniquely determined by $X$ and $g$.
However, the associated triple $(W,\Delta,\alpha)$ of $Y$ is uniquely determined by the associated triple of $X$,
and is computed as in Sec.~\ref{MMPtriples}.
Moreover, if $(W,\Delta,\alpha)$ is terminal, then $Y$ is uniquely determined up to Morita equivalence in $k(X)$
by Cor.~\ref{terminalMoritaequiv}.
\end{construction}

\begin{corollary} \cite[Sec.~4]{CI}
Let $X$ be a terminal order over a projective surface.
Then there exists a Morita equivalence $X \simeq_M X'$ in $k(X)$ and a birational morphism $f \colon X' \rightarrow Y$
to a terminal order $Y$ such that one of the following holds.
\begin{enumerate}
\item $K_Y$ is nef.
\item $Z(Y)$ is the total space of a $\P^1$-bundle over a curve and $K_Y \cdot f < 0$ where $f$ is a fibre.
\item $Z(Y)$ is isomorphic to $\P^2$ and $-K_Y$ is ample.
\end{enumerate}
Moreover, in case (1) $Y$ is uniquely determined up to Morita equivalence in $k(Y)$.
\end{corollary}
\begin{pf}
Let $(Z,\Delta,\alpha)$ be the triple associated to $X$, and $g \colon Z \rightarrow W$ the birational morphism of Thm.~\ref{minmodeltriple}.
Let $\O_Y$ be the reflexive hull of $g_*\O_X$.
Then $Y$ is a terminal order over $W$ satisfying condition (1),(2) or (3) above.
In case (1) $W$ is uniquely determined by Thm.~\ref{minmodeltriple}, and so $Y$ is uniquely determined up to Morita equivalence in $k(Y)$ by
Cor.~\ref{terminalMoritaequiv}.

Let $X' \rightarrow Y$ be a birational morphism obtained from $Y$ and $g \colon Z \rightarrow W$ by the construction~\ref{blowupconstruction}.
Then $X$ and $X'$ are Morita equivalent in $k(X)$ by Cor.~\ref{terminalMoritaequiv}.
\end{pf}

\begin{remark}
Canonical singularities of varieties are the singularities that arise on canonical models of varieties of general type.
We note that the analogous statement is true for orders over surfaces.
Indeed, let $X$ be a terminal order over a projective surface $Z$ with ramification divisor $\Delta$.
If the linear system $|N(K_Z+\Delta)|$ defines a birational map for some $N > 0$ 
(that is, the pair $(Z,\Delta)$ is of \emph{log general type}), we define the \emph{canonical model} $\oX$ of $X$ as follows.
Let
$$\oZ = \Proj \oplus_{n \ge 0} H^0(Z,nK_Z+ \lfloor n \Delta \rfloor).$$ 
(Equivalently, $\oZ$ is the stable image of the map defined by $|N(K_Z+\Delta)|$ for $N$ sufficiently large and divisible.)
Then there is a birational morphism $f \colon Z \rightarrow \oZ$, and we define $\O_{\oX}$ as the reflexive hull of $f_*\O_X$.
One shows that $K_{\oX}$ is ample and the order $\oX$ has canonical singularities, as in the commutative case \cite[p.~355]{R}. 
We do not use this result in this paper.
\end{remark}

\begin{definition}
Let $X$ be a normal order over a surface.
A \emph{resolution} is a proper birational morphism $f \colon Y \rightarrow X$ from a terminal order $Y$.
A resolution $f \colon Y \rightarrow X$ is \emph{minimal} if $K_Y$ is $f_Z$-nef. 
\end{definition} 

\begin{corollary}\label{minresorder}
Let $X$ be a normal order over a surface.
Then $X$ has a minimal resolution $f \colon Y \rightarrow X$ and $Y$ is uniquely determined up to Morita equivalence in $k(X)$.
\end{corollary}
\begin{pf}
Existence follows from Thm.~\ref{minrestriple} and the construction \ref{blowupconstruction}.
Uniqueness is an instance of Cor.~\ref{terminalMoritaequiv}.
\end{pf}

\end{section}

\begin{section}{Artin covers}  \label{frt}
To describe the complete local structure of canonical orders, 
we first classify the possible ramification data and then construct orders with the correct ramification. 
In \cite{A86}, Artin developed a method for constructing normal orders with given ramification data. 
We review his construction here in the format we use to construct canonical orders. 

Let $R$ be a two-dimensional normal noetherian domain with field of fractions $K$ and $A$ a normal $R$-order. 
Let $L/K$ be a Galois field extension with Galois group $G$ and $S$ the integral closure of $R$ in $L$.
 
\begin{theorem} \cite[p.~199, Thm.~2.15]{A86} \label{artincover} 
Suppose that at each curve $C_i \subset \Spec R$, 
the ramification index $r_i$ of $S/R$ at $C_i$ divides the ramification index $e_i$ of $A/R$ at $C_i$. 
Then there is a canonically defined normal $S$-order $B$ satisfying the following properties. 
\begin{enumerate}
 \item The ramification indices of $B$ are $e_i/r_i$ above $C_i$.
 \item $B$ contains $S \otimes_R A$, the $G$-action on $S \otimes_R A$ extends to $B$, and $B^G = A$. 
 \item The natural map $B*G \rtar \End_A B$ is an isomorphism. 
\end{enumerate}
In particular, if $S=k[[u,v]]$ and $e_i = r_i$ for all $i$ then there exists an action of the group $G$ on the full matrix algebra 
$B=k[[u,v]]^{m \times m}$ such that $A = B^G$ and the group action restricts to the Galois action on $S$.
\end{theorem}

We  call the order $B$ the {\it Artin cover} of $A$ with respect to the cover $\Spec S \rtar \Spec R$. 

We recall the notion of reflexive Morita equivalence. 
For a normal $R$-order $A$, let $\refl A$ denote the category of reflexive $A$-modules.
(Here we say an $A$-module is \emph{reflexive} if it is reflexive as an $R$-module.) 
We say two normal orders $A,A'$ are \emph{reflexive Morita equivalent} if there is an equivalence of categories $\refl A \simeq \refl A'$. 
There is a complete analogue of Morita theory in this setting, see \cite[Sec.~1.2]{RVdB} for details. 
We only need the following result.

\begin{lemma}  \label{lsameram}  
Let $A$ be a normal $R$-order and $N$ a reflexive $A$-module such that for every height one prime $P \subset R$, 
\begin{enumerate}
\item $N \otimes_A A_P$ is a generator in the category of $A_P$-modules.
\item $N \otimes_A A_P^{sh}$ is a free $A_P^{\sh}$-module, where $A_P^{\sh}$ is the strict henselisation of $A_P$.
\end{enumerate} 
Then $A':= \End_A N$ is a normal order with the same ramification data as $A$. It is reflexive Morita equivalent to $A$ via the functor $\Hom_A(N,-): \refl A \rtar \refl A'$.
\end{lemma}
\begin{pf}
\'Etale locally at $P$, $A'=\End_A N$ is a full matrix algebra over $A$. Hence $A'$ is normal and has the same ramification data as $A$. The functor
$\Hom_A(N,-) \colon \refl A \rightarrow \refl A'$ is an equivalence by \cite[Prop.~1.2(b)]{RVdB}. 
\end{pf}

We record here some facts about the Artin cover. 

\begin{proposition}  \label{sameram} 
Let $A$ be a normal order and $B$ an Artin cover of $A$. 
Then $A$ and $\End_A B$ are reflexive Morita equivalent normal orders and have the same ramification data. 
\end{proposition}
\begin{pf}
We use Lem.~\ref{lsameram} (note that $B$ is a reflexive $A$-module since it is normal). 
We first localise at a height one prime $P$ in the Zariski topology. 
The order $A = B^G$ is a direct summand of $B$ so certainly $B$ is a generator for $A$ locally at $P$. 
We now further localise in the \'etale topology at $P$. 
We show that $B$ is a free $A$-module. 
Artin gives an \'etale local description of $B$ as follows \cite[p.~202]{A86}. 
Let $r$ be the ramification index of $S/R$ and $e$ the ramification index of $A/R$. 
Then $S=R[s]/(s^r-t)$ where $t \in R$ is a local parameter and
$$ A \otimes_R S= 
\begin{pmatrix}
S      &  \dotsb  & \dotsb  &  S  \\
(s^r)  &  \ddots  &         &  \vdots \\
\vdots &  \ddots  & \ddots  &  \vdots \\
(s^r)  &  \dotsb  & (s^r)   &  S
\end{pmatrix}^{f \times f}$$
We may assume $f=1$.
Let $e_{ij}$ be the standard basis for $S^{e \times e}$ and write $l=\frac{e}{r}$. 
Then 
$$ B = \bigoplus_{i,j} (s^{w_{ij}})e_{ij}$$ 
where the exponents $w_{ij}$ are constant on diagonals 
and are given by the formula
$$ w_{ij} = \left\lceil \frac{i-j}{l} \right\rceil.$$
It is easy to check that $B$ is freely generated as an $A$-module by the 
matrices 
$$\begin{pmatrix}
       &       &         &s^{-c}  &        &   \\
       &       &         &        & \ddots &   \\
       &       &         &        &        & s^{-c} \\
s^{r-c}&       &         &        &        &        \\
       & \ddots&         &        &        &        \\
       &       & s^{r-c} &        &        &
\end{pmatrix}$$
for $0 \le c < r$, where the $s^{-c}$ occur in the $(i,j)$ entries with 
$j-i = cl$ and the $s^{r-c}$ occur in the $(i,j)$ entries with $j-i = cl - e$. 
\end{pf}

\end{section}

\begin{section}{Quotient singularities and skew group rings} \label{sinvariant}
Throughout this section all orders have centre a complete local ring.

Let $S = k[[u,v]]$ and $B = S^{m \times m}$. 
Let $G$ be a finite group acting faithfully on $S$.
Suppose given an extension of the action of $G$ on the centre $S$ of $B$ to an action on $B$.
Then $B^G$ is a reflexive $S^G$-order.
A {\it quotient singularity} is a normal order of this form.  
(We give the explicit condition on the group action for the invariant ring to be normal in Thm.~\ref{normal}.)
In Sec.~\ref{slt}, we show that all log terminal orders have an Artin cover with smooth centre, 
so they are quotient singularities by Thm.~\ref{artincover}. 
Thus, to classify them, we need only compute the possible actions of $G$ on $B$. 
In this section, we first review Artin's approach in \cite{A86} for describing these group actions. 
This is then used to study the skew group ring $B*G$, which is reflexive Morita equivalent to $A=B^G$ when $A \subset B$ is an Artin cover
by Thm.~\ref{artincover}(3) and Prop.~\ref{sameram}. 
These results allow us to interpret log terminal orders as coordinate rings of stacks and to compute their Auslander--Reiten (AR) quivers. 
This material is essentially in  \cite[Sec.~5]{A86} and \cite[Sec.~V.3]{LVV} but, in both cases, the treatment is not quite in the form we want. 

\begin{lemma}\cite[Sec.~4]{A86} \label{projectiverep}
Let $G$ be a finite group acting faithfully on $S=k[[u,v]]$.
Suppose given an extension of the action of $G$ on $S$ to an action on $B=S^{m \times m}=\End_S(S^m)$.
Write $A=B^G$, $R=S^G$, and let $K$ denote the fraction field of $R$.
\begin{enumerate}
\item After a change of coordinates on $S$, the action of $G$ on $S$ is linear, i.e., given by an inclusion $G \subset \GL_2(ku+kv)$. 
\item After a change of basis of $S^m$, the action of $G$ on $B$ is equal to the tensor product of the action on $S$ 
and an action on $k^{m \times m}$ given by a projective representation $\beta \colon G \rtar \PGL_m$.
\item The projective representation $\beta$ lifts to a representation $b \colon G' \rightarrow \GL_m$
of a central extension $G'$ of $G$ by a cyclic group of order $d$, where $d$ is the order of the Brauer class of $A_K$.
\end{enumerate}
\end{lemma}
\begin{pf}
(1) is well-known, see e.g. \cite[p.~97]{Ca}. (2) is \cite[Lem.~4.19]{A86}.
Treating $\PGL_m$ as a trivial $G$-module, we can interpret $\beta$ as an element of the cohomology group $H^1(G, PGL_m)$.
Hence, from the exact sequence 
\begin{equation} 
1 \rtar k^{\times} \rtar \GL_m \rtar \PGL_m \rtar 1  
\end{equation}
we obtain a class $d\beta \in H^2(G,k^{\times})$.
Then $d\b$ is determined by the Brauer class of $A_K$ and has the same order $d$ by \cite[Prop.~4.9]{A86}.
By Kummer theory, $d\beta$ is induced from an element $\gamma \in H^2(G,\mu_d)$. 
The class $\gamma$ defines a central extension
$$ 
1 \rightarrow \mu_d \rightarrow G' \rightarrow G \rightarrow 1
$$
such that $\beta$ lifts to a representation $b \colon G' \rightarrow \GL_m$.
\end{pf}

\begin{remark} \label{centralextensionremark}
We often lift the projective representation $\b$ to a representation $b$ of a central 
extension $G'$ of $G$ by $\mu_d$ with $d$ some multiple of the minimal value determined above.
Then the image of $G'$ in $\GL_m$ contains the $d$th roots of unity.
This is used to simplify our calculations in Sec.~\ref{exp}. 
\end{remark}


%
%
%
%
%
%
%

\vspace{2mm}

We now consider the skew group ring $B*G$. The following basic facts are known to experts.

\begin{lemma} \label{lmoritaskew}
Let $\eta$ be a generator of the central subgroup $\mu_d \subset G'$.
\begin{enumerate}
 \item There is a Morita equivalence between $B*G'$ and $S*G'$ given 
by the Morita bimodule $P:=(S*G')^m$. The right action of $S*G'$ is given 
by coordinatewise multiplication and the left action of $B$ is given by 
matrix multiplication. 
Left multiplication by $g \in G'$ is given by coordinatewise left 
multiplication by $g$ followed by left multiplication by the matrix $b_{g}$. 
 \item In $B*G'$ we have a central idempotent $f:=\frac{1}{d}(1+\eta + \ldots + \eta^{d-1})$. 
Hence $B*G' = f(B*G') \times (1-f)(B*G')$ and we have natural isomorphisms of algebras $f(B*G') \simeq B*G'/(1-\eta) \simeq B*G$. 
Similarly, if $\z$ is a $d$-th root of unity then $e:=\frac{1}{d}(1+\z\eta + \ldots + \z^{d-1}\eta^{d-1})$ is a central idempotent 
in $S*G'$ and we have a natural isomorphism $e(S*G') \simeq S*G'/(1-\z\eta)$. 
 \item The Morita equivalence in (1) restricts to a Morita equivalence between $f(B*G')$ and $e(S*G')$ via the Morita bimodule $fP = Pe$, where $\zeta$ is defined by $b_{\eta}=\zeta I$. 
\end{enumerate}
\end{lemma}
\begin{pf}
For $(1)$ we first prove the case where the action $b$ is trivial. In this case we denote the skew group ring by $BG'$. 
Then $BG' = (S*G')^{m \times m}$ so the Morita equivalence is clear. 
When the action is arbitrary we have an isomorphism $\phi:B*G' \xrightarrow{\sim} BG'$ which is the identity on $B$ and 
maps $g$ to $b_{g}g$. 
Hence the Morita equivalence for $BG'$ induces the given one for $B*G'$. 

Note that $\eta$ is central in both $B*G'$ and $S*G'$ so $(2)$ and $(3)$ follow from explicit computations or standard Morita theory. 
\end{pf}

We have the following different interpretations of quotient singularities which are due primarily to Le Bruyn, Reiten, Van den Bergh and Van Oystaeyen.

\begin{theorem} (\cite[Prop.~V.3.2]{LVV}, \cite[Thm.~5.6]{RVdB}) \label{twistedskewgroupring}
\label{skewgroup}
Let $A$ be a quotient singularity and $\a$ the Brauer class of $k(A)$. 
Then there are 
\begin{enumerate}
\item a finite group $G \subset \GL(ku+kv)$, 
\item a central extension $G'$ of $G$ by a cyclic group, and
\item a central irreducible idempotent $e$ of $kG'$
\end{enumerate} 
such that $A$ is reflexive Morita equivalent to the algebra $e(S* G')$ of global dimension $2$. 
Equivalently, $A$ is reflexive Morita equivalent to the twisted skew group ring $S*_c G$ where $c \in H^2(G,k^{\times})$ is the class corresponding to $\a$.
\end{theorem}

\begin{pf}
The order $A$ has an Artin cover $A \subset B$ where $B=S^{m \times m}$ by \cite[Pf. of Thm.~5.6]{RVdB} (or the proof of our Thm.~5.4).  
The orders $A$ and $B*G$ are reflexive Morita equivalent by Prop.~\ref{sameram} and Thm.~\ref{artincover}. 
We construct the central extension $G'$ as in Lem.~\ref{projectiverep}. 
The orders $B*G$ and $e(S* G')$ are Morita equivalent by Lem.~\ref{lmoritaskew}.  
The last statement is \cite[Prop.~V.3.2]{LVV}, and can also be proved directly by showing that $e(S* G') \simeq S*_c G$.
\end{pf}

%

Thm.~\ref{twistedskewgroupring} gives a way to compute AR quivers of quotient singularities from McKay quivers. 
Here, for $W$ a representation of a group $G$, the {\it McKay quiver} of $G$ acting on $W$ is the 
quiver with vertices labelled by the irreducible representations $V_i$ of $G$ and 
the number of arrows from $V_i$ to $V_j$ given by $\dim_k \Hom_G(V_i, W \otimes_k V_j)$.
For more background on McKay quivers see \cite{Ausmck}.


\begin{proposition} \label{arquivers}
Let $A$ be a quotient singularity, $G'$ the finite group and $b \colon G' \rightarrow \GL_m$ the representation 
determined in Thm.~\ref{twistedskewgroupring}. Let $\eta$ be a generator of the central subgroup $\mu_d \subset G'$.
Write $b_{\eta}=\zeta I$ where $\zeta$ is a $d$th root of unity.
Then the AR quiver of $A$ is the component of the McKay quiver of $G'$ acting on $V=ku+kv$
where $\mu_d$ acts via the character $\eta \mapsto \zeta^{-1}$.
\end{proposition}
\begin{pf} 
We use the notation of Lem.~\ref{lmoritaskew}. 
Auslander's proof of the McKay correspondence \cite[Sec.~2]{Aus} shows that there is a correspondence between indecomposable projective modules 
over $S*G'$ and irreducible $G'$-modules which identifies the AR quiver of $S*G'$ and the McKay quiver of $G'$ acting on $V$. 
The correspondence is given by tensoring $G'$-modules with $S$ and taking fibres of $S*G'$-projectives at the closed point. 
Since $\mu_d \subset G'$ is central and acts trivially on $V$, the McKay quiver splits up as a union of components 
where $\mu_d$ acts by a given character. 
The component of the McKay quiver where $\eta$ acts by $\z^{-1}$ thus corresponds to the AR-quiver of $S*G'/(1-\z\eta)$. 
The algebra $S*G'/(1- \zeta\eta)$ is Morita equivalent to $B*G$ by Lem.~\ref{lmoritaskew}. 
Finally $B*G$ and $A=B^G$ are reflexive Morita equivalent by Thm.~\ref{artincover}(3) and Prop.~\ref{sameram} so they have the same AR-quiver.
\end{pf}

\end{section}

\begin{section}{Log terminal orders} \label{slt}

Let $(P \in Z)$ be the spectrum of a complete local normal domain of dimension $2$ and $X$ a normal $Z$-order. 
We say $X$ has {\it finite representation type} if the number of isomorphism types of indecomposable reflexive modules is finite.
By \cite[Ch.~V]{LVV}, the normal orders with finite representation type are the quotient singularities. 
In this section, we show that these orders can also be characterized geometrically as the log terminal orders.

We recall the construction of a noncommutative version of the index one cover 
(\cite[Sec.~V.2]{LVV}, see \cite[Def.~5.19]{KM} for the commutative case). 
Let $X$ be a normal order over a surface germ $(P \in Z)$ as above. 
For $L$ an $\O_X$-bimodule and $N \in \N$, define $L^{[N]}=(L^{\otimes N})^{**}$, 
the reflexive hull of the $N$th tensor power of $L$ over $\O_X$. 
By Prop.~\ref{candiv} we have a natural isomorphism of $\O_X$-bimodules 
$$
\w_X^{[m]} \simeq \O_X \otimes \O_Z(m(K_Z+\Delta)),
$$
where $m$ is the degree of $X$. 
Assume that the $\Q$-divisor $K_Z+\Delta$ is $\Q$-Cartier. 
Then 
$\w_X^{[N]} \simeq \O_X$ as bimodules for some $N \in \N$. 
Let $n$ be the least such $N$, the \emph{order} of $\w_X$.
We define the {\it index one cover} $\tilde{X}$ of $X$ by
$$ 
\O_{\tilde{X}} := \O_X \oplus \w_X \oplus \w_X^{[2]} \oplus \cdots \w_X^{[n-1]},
$$
where the multiplication is given by fixing an isomorphism $\w_X^{[n]} \cong \O_X$.
 
A proof of the following result was sketched in \cite[Cor.~V.2.3, Prop.~V.2.4]{LVV}.
Here we include a direct proof.

\begin{proposition} \label{lcancodim1}
Let $X$ be a normal order over a surface germ $(P \in Z)$ such that $K_Z+\Delta$ is $\Q$-Cartier. 
Let $\tilde{X}$ be the index one cover. Then $\O_{\tilde{X}}$ is Azumaya in 
codimension one and the centre $\tilde{Z}$ of $\tilde{X}$ has the 
same ramification indices over $Z$ as $\O_X$. In particular, 
$\O_{\tilde{X}}$ is an Artin cover of $\O_X$. Moreover,
$\w_{\tilde{X}}$ is isomorphic to $\O_{\tilde{X}}$ as a bimodule 
and $\tilde{Z}$ is Gorenstein.
\end{proposition}
\begin{pf}
We may work \'etale locally in codimension one.
Let $e$ be the ramification index of $\O_X$ at some codimension 
one prime. 
Let $R$ be the centre of $\O_X$, a discrete valuation ring, and $t$  
a uniformising parameter. We may assume $\O_X$ has the standard form 

\begin{equation}   \label{hinv}
\O_X = 
\begin{pmatrix} 
R      &   \hdots  & \hdots & R  \\
(t)    &   R       &        &  \vdots       \\
\vdots &  \ddots   & \ddots &              \\ 
(t)    &  \hdots   &  (t)   &   R
\end{pmatrix}^{f \times f} \subset R^{ef \times ef}.
\end{equation}
where the displayed matrix is $e \times e$.
Locally, the order of $\w_X$ equals $e$. 
To simplify notation we assume $f=1$. 

Let $u$ be an $e$-th root of 
$t$ and $T:= R[u]$. We show that 
$\O_{\tilde{X}}$ is $\mu_e$-equivariantly isomorphic to the matrix algebra   
\[\O_{\bar{X}} = 
\begin{pmatrix}
T      & u^{-1}T   & u^{-2}T  & \hdots  \\
uT     &  T        & u^{-1}T  & \hdots  \\
u^2T   & uT        & T          & \hdots  \\
\vdots & \vdots    &  \vdots    & \ddots
\end{pmatrix} \simeq T^{e\times e}\]
where $\mu_{e}$ acts on $\O_{\bar{X}}$ via the the natural Galois action on $T$ over $R$.

The bimodule $\w_X$ is generated as a left and right $\O_X$-module by the element 
\[ x = 
\begin{pmatrix}
0 & t^{-1}  &        &   \\
  & \ddots  & \ddots &  \\
  &         & \ddots & t^{-1}  \\
1 &         &        &  0
\end{pmatrix}.\]
Since we are working \'etale locally, we may assume that the 
relation in the canonical cover is given by 
$x^e = 1$.  Then there is a $\mu_e$-equivariant isomorphism of $\O_X$-algebras
$\O_{\tilde{X}} \xrightarrow{\sim} \O_{\bar{X}}$ 
given by 
\[ x \mapsto y = 
\begin{pmatrix}
0 & u^{-1}  &        &   \\
  & \ddots  & \ddots &  \\
  &         & \ddots & u^{-1}  \\
u^{e-1} &         &        &  0
\end{pmatrix}.\]
Indeed, $y^e = 1$, $y$ commutes with $\O_X$ the same way $x$ does, and 
the powers of $y$ generate the $\mu_e$-eigenspaces of $\O_{\bar{X}}$. 

We now show $\w_{\tilde{X}}$ is isomorphic to $\O_X$ as a bimodule.
We have $\w_{\tilde{X}}=\Hom_{\tilde Z}(\O_{\tilde X},\w_{\tilde Z})=\Hom_Z(\O_{\tilde X},\w_Z)$,
and $\O_{\tilde{X}}=\O_X \oplus \w_X \oplus \cdots \oplus \w^{[n-1]}_X$ as a $\O_X$-bimodule.
So $\w_{\tilde{X}}$ contains the direct summand $\Hom_Z(\w_X,\w_Z)=\Hom_Z(\Hom_Z(\O_X,\w_Z),\w_Z)$,
which is canonically isomorphic to $\O_X$ as a bimodule. Let $s$ be the generator corresponding to $1 \in \O_X$.
One checks that $s$ is central in the $\O_{\tilde{X}}$-bimodule $\w_{\tilde{X}}$. Indeed, $\O_{\tilde{X}}$ is generated over 
$\O_X$ by $\w_X \subset \O_{\tilde{X}}$, and, for $a \in \w_X$, we have $as=sa=i(a)$, where $i : \w_X \rightarrow \w_{\tilde{X}}$
is the natural inclusion. Thus $\w_{\tilde{X}}$ is isomorphic to $\O_{\tilde{X}}$, generated by $s$.
Finally, passing to the centres of these bimodules, we deduce $\w_{\tilde{Z}} \cong \O_{\tilde Z}$, i.e., 
$\tilde{Z}$ is Gorenstein.
\end{pf}


Let $Z$ be a normal variety and $\Delta$ an effective $\Q$-divisor on $Z$.
If $f \colon W \rightarrow Z$ is a resolution of $Z$ we have an equality 
$K_W+\Delta'=f^*(K_Z+\Delta)+\sum a_i E_i$ where $\Delta'$ is the strict transform of $\Delta$, the $E_i$ are the exceptional curves, and $a_i \in \Q$.
The \emph{discrepancy} $\discrep(Z,\Delta)$ of $(Z,\Delta)$ is the infimum of the coefficients $a_i$ over all 
exceptional divisors of all resolutions of $Z$.
We say $(Z,\Delta)$ is \emph{(Kawamata) log terminal} if $\discrep(Z,\Delta) > -1$ and the coefficients of $\Delta$
are less than $1$.

\begin{proposition} \label{pfrtislt}
Let $X$ be a quotient singularity and $(Z,\Delta,\alpha)$ the associated triple. Then
\begin{enumerate}
\item $(Z,\Delta)$ is log terminal. 
\item $\discrep(Z,\Delta,\alpha) \geq \discrep(Z,\Delta)$.
\end{enumerate}
In particular, $X$ is log terminal.
\end{proposition}
\begin{pf}
Let $\tilde{X}$ be the index one cover of $X$ and $\tilde{Z}$ the centre of $\tilde{X}$. 
Then $\tilde{X}$ has finite representation type by \cite[Prop.~2.12]{A86} since $\tilde{X}$ is an Artin cover of the
quotient singularity $X$ by Prop.~\ref{lcancodim1}. 
It follows that $\tilde{Z}$ has finite representation type and so is a quotient singularity.
Equivalently, $\tilde{Z}$ is log terminal \cite[Prop.~4.18]{KM}.
Consider the induced cover $\pi \colon \tilde{Z} \rtar Z$. 
The ramification indices of this cover are the same as those of $X$ by Prop.~\ref{lcancodim1}, 
so $K_{\tilde{Z}} = \pi^*(K_Z + \D)$ by the Riemann--Hurwitz formula. 
Hence $(Z,\D)$ is log terminal by \cite[Prop.~5.20]{KM}. 

There exists a single resolution  $f \colon W \rightarrow Z$ of $(Z,\Delta,\alpha)$ which can be used to compute $\discrep(Z,\Delta,\alpha)$ and 
$\discrep(Z,\Delta)$ by Prop.~\ref{computediscrep} and \cite[Cor.~2.32(2)]{KM}.
Write 
$$ K_{W} + \Delta_W = f^*(K_Z + \Delta) + \sum a_i E_i,$$
$$ K_{W} + \Delta'= f^*(K_Z + \Delta) + \sum b_i E_i.$$
where the $E_i$ are the exceptional curves.
If $e_i$ is the ramification index of $(W,\Delta,\alpha)$ along $E_i$ 
then $a_i =( 1 - \frac{1}{e_i}) + b_i$.
To prove (2), it suffices to show that $e_ia_i \ge b_i$ for all $i$, equivalently, 
$(b_i+1)(e_i - 1) \geq 0$. But $(Z,\D)$ is log terminal so $b_i > -1$.
\end{pf}

\begin{remark}
Prop.~\ref{pfrtislt}(2) has the following stacky interpretation. 
Consider an arbitrary quotient singularity which, by Thm.~\ref{skewgroup}, is reflexive Morita equivalent 
to and has the same ramification data as $\O_X:=e(k[[u,v]]*G')$, where $G'$ is a central extension of some finite 
subgroup $G \subset GL_2$ by $\mu_d$ and $e$ is a central idempotent. 
Let $\mathcal{S}'$ denote the Deligne-Mumford stack $[\Spec k[[u,v]]/G']$.
Then there is a fully faithful functor from the category $\Mod X$ of $\O_X$-modules to the category of quasi-coherent sheaves on $\mathcal{S}'$ 
with image the sheaves where the stabilizer $\mu_d$ acts by a fixed character determined by $e$. 
The stack $\mathcal{S}'$ is a $\mu_d$-gerbe over the orbifold $\mathcal{S}=[\Spec k[[u,v]]/G]$, 
and the noncommutative coordinate ring of $\mathcal{S}$, as defined in \cite{CI1}, is the order $\O_Y = k[[u,v]] * G$. 
In these terms, Prop.~\ref{pfrtislt}(2) states that the discrepancy of $\mathcal{S}'$ is greater than or equal to that of $\mathcal{S}$. 
\end{remark}



\begin{theorem} \label{tltisfrt}
Let $X$ be a normal order over the spectrum $(P \in Z)$ of a complete local normal domain of dimension $2$.
The following are equivalent. 
\begin{enumerate}
\item $X$ has finite representation type.
\item $X$ is a quotient singularity.
\item $X$ is log terminal.
\end{enumerate}
\end{theorem}
\begin{pf}
$(1)$ and $(2)$ are equivalent by \cite[Ch.~V]{LVV} and $(2)$ implies $(3)$ by Prop.~\ref{pfrtislt}.
Let $X$ be a log terminal order with centre $Z$ and ramification divisor $\Delta$. 
We show that $X$ is a quotient singularity. 
Let $\tilde{X}$ be the index one cover of $X$, $\tilde{Z}$ the centre of $\tilde{X}$, and $\pi \colon \tilde{Z} \rightarrow Z$ the induced map. 
Then $\tilde{Z}$ and $X$ have the same ramification indices over $Z$ by Prop.~\ref{lcancodim1}. 
Hence $K_{\tilde{Z}} = \pi^*(K_Z + \D)$ by the Riemann-Hurwitz formula. 
Since $X$ is log terminal, the pair $(Z,\D)$ is log terminal by \cite[Prop.~3.15]{CI}.  
Hence $\tilde{Z}$ is also log terminal by \cite[Prop.~5.20(3)]{KM}, i.e., $\tilde{Z}$ is a quotient singularity.
Let $q \colon U \rightarrow \tilde{Z}$ be the smooth cover of $\tilde{Z}$ and $B$ the reflexive hull of the pullback of $\O_{\tilde X}$. 
Then $B$ is  Azumaya in codimension $1$ because $\O_{\tilde{X}}$ is so and $q$ is \'etale in codimension $1$.
Thus $B$ is Azumaya because $U$ is smooth, and so $B \simeq \O_U^{m \times m}$.
The composite map $U \rightarrow Z$ is Galois with some group $G$ (cf. \cite[Lem.~4.2]{A86}), and the $G$-action extends to $B$ so that $\O_X=B^G$. 
Hence $X$ is a quotient singularity as required.
\end{pf}

In Sec.~\ref{exp} we explicitly describe all canonical orders as invariant rings 
$(k[[u,v]]^{m \times m})^G$.

\end{section}

\begin{section}{Minimal resolutions of canonical orders}\label{sminres}

In this section, we classify the ramification data of canonical singularities of orders. 
Our approach is to first determine their minimal resolutions.

\begin{proposition} \label{tminres}
Let $X$ be a canonical order and $f \colon Y \rightarrow X$ the minimal resolution. Then $K_Y=f_Z^*K_X$.
\end{proposition}
\begin{pf}
Write $K_Y=f_Z^*K_X+\sum a_iE_i$ where the $E_i$ are the exceptional curves and $a_i \in \Q$.
Then $a_i \ge 0$ because $X$ is canonical and $a_i \le 0$ because $f$ is minimal by Rem.~\ref{minresremark}.
\end{pf}

\begin{proposition} \label{4cases}
Let $X$ be a canonical order and $E$ an exceptional curve of the minimal resolution $f \colon Y \rightarrow X$. 
Then $E$ is a smooth rational curve. Let $D = \Supp \D$ be the discriminant of $Y$. 
Then one of the following holds in an \'etale neighbourhood of $E$.
Here $U$ and $V$  denote ramification curves that meet $E$ transversely at two distinct points.
\begin{eqnarray*}
\renewcommand{\arraystretch}{1.5}
\begin{array}{|l|l|l|}
\hline 
					& E^2=-1					& E^2=-2 			\\
\hline
		E \not\subset D 	& \D=\frac{1}{2} D, \, E \cdot D=2 		& E \cap D = \emptyset \\
\hline
		E \subset D     	& \D=(1-\frac{1}{e})E+(1-\frac{1}{2e})(U+V) 	& \D=(1-\frac{1}{e})(E+U+V) \\
\hline
\end{array}
\end{eqnarray*}
\end{proposition}

\begin{pf}  
Let $Z$ denote the centre and $\Delta$ the ramification divisor of $Y$.
We have $E^2<0$ because $E$ is exceptional and $E \cdot K_Y= E \cdot (K_Z + \Delta) = 0$ by Prop.~\ref{tminres}.  

Suppose first that $E \not\subset D$.  
Then $E \cdot \D \geq 0$, so $E \cdot K_Z \leq 0$ and $E \cdot (K_Z+E)<0$.
Hence $E$ is a smooth rational curve with $E^2=-1$ or $E^2=-2$, and $E \cdot \D=1$ or $E \cdot \D=0$ respectively.  
The coefficients of $\D$ are at least $\frac{1}{2}$ so we obtain the cases above.  

Now suppose that $E \subset D$ and write $\D=(1-1/e)E+\D'$ and $D=E+D'$.
Then $$0=E \cdot (K_Z+\D)=E \cdot (K_Z+E) -\frac{1}{e} E^2 + \D' \cdot E > E\cdot(K_Z+E) ,$$
so $E$ is a smooth rational curve. Hence $\D' \cdot E =2 + \frac{1}{e}E^2 <2$.
We proceed as in \cite[Pf. of Thm.~3.10]{CI}, using the local description of the ramification data of terminal orders in Thm.~\ref{tterminal}. 
We have $D' \cdot E \le 3$. If $D' \cdot E = 3$  then  
$$\D' \cdot E= \left(1-\frac{1}{e_1}\right) + \left(1 - \frac{1}{e_2}\right) +\left(1- \frac{1}{e_3}\right) < 2$$
where the $e_i$ are the ramification indices of the ramification curves intersecting $E$. 
The solutions $(e_1,e_2,e_3)$ of the above inequality are the Platonic triples. 
None of the solutions yield ramification indices for a cyclic cover of $E \simeq \P^1$, a contradiction. So $D' \cdot E \le 2$.  
A nontrivial cyclic cover of $E$ ramifies over at least $2$ points, so $D' \cdot E =2$. 
The two ramification curves intersecting $E$ must have ramification indices $me,ne$ for some $m,n \in \N$. 
Then $(1-\frac{1}{ne}) + (1 - \frac{1}{me})= 2 + \frac{1}{e}E^2$, or $-E^2 = \frac{1}{n}+\frac{1}{m}$.
The only solutions are $n=m=1$ and $n=m=2$ for $E^2=-2$ and $-1$ respectively. 
\end{pf}

\begin{corollary} \label{corom2}
Let $E_1,E_2$ be exceptional curves of the minimal 
resolution $Y \rtar X$.  Suppose that $E_1,E_2$ intersect and 
that $E_1^2 = -2$. 
\begin{enumerate}
\item If $E_2^2 = -2$ then $Y$ ramifies along $E_1$ if and 
only if it ramifies along $E_2$. Furthermore, the ramification 
indices on $E_1,E_2$ are the same. 
\item If $E_2^2 = -1$ then $Y$ is unramified along $E_2$. 
\end{enumerate}
\end{corollary}
\begin{pf}
(1) follows from Prop.~\ref{4cases} applied to $E = E_1,E_2$. 
Suppose now that $E_2^2 = -1$.  
If $E_2 \subset D$ then Prop.~\ref{4cases} applied to $E=E_1$ shows that $E_1 \subset D$. 
Let $e_1,e_2$ be the corresponding ramification indices.
Applying Prop.~\ref{4cases} again to $E=E_1,E_2$ in turn yields $e_2 = e_1$ and $e_2 = 2e_1$, a contradiction. 
\end{pf}

\begin{lemma}\label{expcurves}
Let $E=\bigcup E_i$ be a connected union of curves on a smooth surface such that each $E_i$ is a $(-1)$-curve or a $(-2)$-curve, 
and the matrix $(E_i.E_j)_{ij}$ is negative definite.
Then $E$ is either an ADE configuration of $(-2)$-curves, 
or a chain of curves such that one end of the chain is a $(-1)$-curve and the remaining curves are $(-2)$-curves. 
In particular, the exceptional locus of the minimal resolution $Y \rtar X$ of a canonical order is of this form. 
\end{lemma}
\begin{pf}
If the curves are all $(-2)$-curves then it is well known that they must form an ADE configuration.
If $E_1$ is a $(-1)$-curve then we can contract $E_1$ and argue by induction. 
\end{pf}

\begin{theorem}  \label{canramdata}
Let $X$ be a canonical order over a surface germ $(P \in Z)$.
Assume that $X$ is not terminal and let $f \colon Y \rightarrow X$ be the minimal resolution.

The possible configurations of exceptional and ramification curves on $Z(Y)$ are as listed in the appendix.
They fall into types $A_n$, $D_n$, $E_n$, $A_{1,2,\xi}$, $A_{n,\xi}$, $B_n$, $BD_n$, $L_n$, $BL_n$, $DL_n$, where the 
(first) subscript is the number of exceptional curves of the minimal resolution.
Here $n \ge 1$ for $A_{n,\xi}$, $B_n$, and $L_n$, and $n \ge 2$ for $BD_n$ and $DL_n$. 

The types $A_n$, $D_n$, $E_n$ correspond to orders $X$ over Kleinian singularities which are unramified in codimension one. 
The ramification data of $X$ for the remaining types are listed in the table below.
We list the ramification indices $e_C$ of the components $C$ of the discriminant and, for each $C$, 
the ramification index $e_P$ of the cyclic cover $\tilde{C}$ of the normalisation of $C$ over $P$.
If the centre $Z$ is smooth we let $x$,$y$ denote coordinates at $P \in Z$. 
\begin{eqnarray*} \label{canramtab}
\renewcommand{\arraystretch}{1.5}
\begin{array}{|l|l|l|l|l|} \hline
\mbox{type} 		& \mbox{centre} 	& \mbox{discriminant} 		& e_{C}  	& e_P	\\
\hline
A_{1,2,\xi}  		& \mbox{smooth} 	& xy        			& 2e,2e 	& e,e 	\\
\hline
BL_n  			& \mbox{smooth} 	& y^2-x^{2n+1}  		& 2 		& 1 	\\
\hline
B_n  			& \mbox{smooth} 	& (y+x^n)(y-x^n)         	& 2,2 		& 1,1 	\\
\hline
L_n  			& \mbox{smooth} 	& (y+x^{n+1})(y-x^{n+1})        & 2,2 		& 2,2 	\\
\hline
DL_n  			& \mbox{smooth} 	& x(y^2-x^{2n-1}) 		& 2,2 		& 2,2 	\\
\hline
BD_n  			& \mbox{smooth} 	& x(y+x^{n-1})(y-x^{n-1})	& 2,2,2 	& 2,2,1 \\
\hline
A_{n,\xi}	   	& A_n=(xy=z^{n+1})      & z     			& e,e 		& e,e 	\\
\hline
\end{array}
\end{eqnarray*}

\end{theorem}

\begin{remark}
The names of the types are chosen to agree with \cite{HPR}. We also write $C_2:=BD_2$ for the same reason.
\end{remark}

\begin{pf}  
Let $E$ be the exceptional locus of the minimal resolution $f \colon Y \rightarrow X$. 
If all the components of $E$ are $(-2)$-curves then either $E \cap D = \varnothing$ or $E \subset D$ by Cor.~\ref{corom2}(1).
If $E \cap D = \varnothing$ then $E$ is the exceptional locus of a resolution of a Kleinian singularity. 
If $E \subset D$ then each component $E_i$ of $E$ must intersect exactly two other components of $D$ by Prop.~\ref{4cases}. 
This gives case $A_{n,\xi}$.

Now suppose $E_1$ is a $(-1)$-curve.
Then $E$ is a chain of curves with first curve $E_1$ and the remaining curves being $(-2)$-curves by Lem.~\ref{expcurves}. 
If $E_1 \subset D$ then $E_1$ is the unique exceptional curve by Cor.~\ref{corom2}(2). 
This gives case $A_{1,2,\xi}$ by Prop.~\ref{4cases}. If $E_1 \not\subset D$ 
then the ramification index of each component of $D$ meeting $E$ equals $2$.
The union of the remaining components of $E$ must be either contained in $D$ or disjoint from $D$ by Cor.~\ref{corom2}(1).  
In the first case, all ramification indices are equal to $2$. 
Let $E_2$ denote the exceptional curve in the chain adjacent to $E_1$.  
There is another non-exceptional ramification curve $U$ intersecting $E_1$ transversely because $E_1 \cdot D=2$. 
If $U$, $E_1$, and $E_2$ are concurrent then we are in case $DL_n$, otherwise we are in case $BD_n$. 
Finally suppose $E - E_1$ is disjoint from $D$. 
Since $E_1 \cdot D=2 $ and $D$ has normal crossings, 
$E_1$ meets $D$ either in two distinct points, at a node, or tangentially at a smooth point.
This gives cases $B_n$, $L_n$, and $BL_n$ respectively.

Conversely, if an order $X$ has ramification data as listed in the table, we compute that 
the minimal resolution $f \colon Y \rightarrow X$ is as described above, cf. Rem.~\ref{blowupram}.
Then $K_Y=f_Z^*K_X$, so $X$ is canonical
\end{pf}

\end{section}

\begin{section}{The Gorenstein property}\label{gorprop}

We say that an order $X$ over a surface $Z$ is \emph{Gorenstein} if the dualising sheaf $\omega_X :=\cHom_Z(\O_X,\omega_Z)$ 
is locally isomorphic to $\O_X$ as a left and right $\O_X$-module. 
Recall that a canonical order $A$ has the form $B^G$ where $B = k[[u,v]]^{m \times m}$. 
In this section, we first determine explicitly the condition on the group action for the invariant ring $B^G$ to be normal 
(or, equivalently, Gorenstein in codimension one). We then describe the Gorenstein condition in terms of the group action. 
In Sec.~\ref{exp} we verify this condition case-by-case for canonical orders.



We use the notation of Sec.~\ref{sinvariant}. 
Let $S=k[[u,v]]$ and suppose that $G \subset GL(ku+kv)$ acts on $B = S \otimes k^{m \times m}$ as the tensor product of the natural action on $S$ 
and some action on $k^{m \times m}$. Let $R$ be the centre of $A = B^G$. 
Artin gives the condition on the group action for $A$ to be a maximal order \cite[p.~210, 4.15]{A86}. 
Here we give the condition for $A$ to be normal.  

We may assume that $G$ acts linearly via unitary matrices on $\Spec S$. 
Let $L \subset \Spec S$ be a line of reflection of $G$, and $C$ its image in $\Spec R$.
Let $\hat{(\cdot)}$ denote completion at $L$ or $C$. 
Let $H$ be the stabiliser of $L$, then $\hat{A} = \hat{B}^H$. 
Write $L=(t=0)$ and $L^{\perp}=(s=0)$. Then $H$ acts diagonally on $k s + k t$. 
The inertial group $I=\{h \in H \ | \ h(s) = s\}$ is cyclic, say $I=\langle \sigma \rangle$, and $H/I$ is also cyclic. 
For the subgroups $G$ we are interested in, $H$ is always a split extension of $I$. 
We therefore assume there is a cyclic subgroup $J=\langle \tau \rangle$ of $H$ such  that $H = I \times J$. 
(We verify that this holds for all the groups $G$ in Sec.~9 by showing that $\det \s$ generates $\det(G)$ so the determinant map gives a splitting of the inclusion $I \subset H$.) 

\begin{theorem} \label{normal}
With notation as above, let the action of $\s,\t$ on $k^{m \times m}$ be given by conjugation by $b_\s,b_\t \in k^{m \times m}$. 
Then $\hat{A}=\hat{B}^H$ is a normal order if and only if all the eigenspaces of  $b_\s$ 
have the same dimension. In this case, the ramification index of $A$ along 
$C$ is the number of distinct eigenvalues of $b_{\s}$.
Moreover, let $\eta$ be the scalar such that $b_\s b_\t = \eta b_\t b_\s$. 
Then $\eta$ is a primitive $e$-th root of unity where $e$ is the index of the 
central simple algebra $k(\hat{A})$. Equivalently, $e$ is the ramification index
of the cyclic cover of the normalisation of $C$ determined by $A$ over the closed point.
\end{theorem}

\begin{pf}
For the proof of normality, first note that by \cite[Thm.~7.8.8]{McR}, $\hat{B}^H$ is hereditary. 
So we need only check that the eigenspace condition in the statement is equivalent to the condition that $\hat{B}^H/ \mbox{rad}\ \hat{B}^H$ is a product of matrix algebras of the same size. 

Note that $\s,\t$ commute so $b_\s,b_\t$ must skew commute by some scalar $\eta$ which is a root of unity of some order $e$.
The skew commutation relation shows that $\eta^{|I|} = 1 = \eta^{|J|}$.
Write $|I|=pe$ and $|J|=qe$. 
Let $\z$ be a primitive $pe$-th root of unity and write $\eta = \z^{pl}$. 

We view $b_\s, b_\t$ as endomorphisms of $V = k^m$. 
For convenience, we scale $b_{\s},b_{\t}$ so that $b_{\s}^{pe} = 1 = b_{\t}^{qe}$. 
We consider the eigenspace decomposition $V = \oplus_{i \in \Z/pe\Z}\, V_i$ with respect to $b_\s$, where $V_i$ is the $\z^i$-eigenspace of $b_\s$. 
Write $d_i = \dim V_i$. 
With respect to an appropriate basis, we have 
\begin{equation}  \label{bs} 
b_\s = 
\begin{pmatrix}
I      &   0     &   0    &   \hdots  &  0 \\
0      & \z I    &   0    &   \hdots  &  0 \\
0      &   0     & \ddots &           &  \vdots \\
\vdots & \vdots  &        &  \ddots   &  0      \\
0      &   0     & \hdots &    0      & \z^{pe-1} I   
\end{pmatrix}
\end{equation}
where the $I$ in column $(j+1)$ denotes the identity matrix of size $d_j$. 

Let $T:= \hat{S}^I = k((s))[[t^{pe}]]$. 
We compute $\hat{B}^H/ \rad \hat{B}^H$ by first examining $\hat{B}^I$. 
Since $I$ is the inertia group and the group action is faithful, we may assume $\s$ maps $t$ to $\z t$. 
Computing invariants with respect to $I$ gives the block matrix form 

\[\hat{B}^I = 
\begin{pmatrix}
T      & t^{pe-1}T & t^{pe-2}T  & \hdots  \\
tT     &  T        & t^{pe-1}T  & \hdots  \\
t^2T   & tT        & T          & \hdots  \\
\vdots & \vdots    &  \vdots    & \ddots
\end{pmatrix} \]
where the block sizes are the same as for $b_{\s}$.
We conjugate by the block matrix 
\[ \xi:= \begin{pmatrix}
1      &  0     &  0  &  \hdots  \\
0      &  t     &  0  &  \hdots  \\
0      &  0     & t^2 &     \\
\vdots & \vdots &     &  \ddots
\end{pmatrix}\] 
to obtain 
\begin{equation}   
\begin{pmatrix} 
T      & t^{pe}T   & \hdots & t^{pe}T  \\
T      &   T       & \ddots &  \vdots       \\
\vdots &           & \ddots & t^{pe}T             \\ 
T      &  \hdots   &        &   T
\end{pmatrix}
\end{equation}

Then 
\[ \hat{B}^I/\rad \hat{B}^I \cong 
   \prod_{i=0}^{pe-1} k((s))^{d_i \times d_i}   .\]

To pass to the residue ring of $\hat{B}^H$ we need the following.
\begin{lemma}   \label{resinv} 
Let $C$ be a ring and suppose a group $J$ acts 
by automorphisms on $C$. Then 
$\rad C^J \supseteq \rad C \cap C^J$. 
\end{lemma}
\begin{pf}
Let $r \in \rad C \cap C^J$. 
We need to show that $1-r$ is invertible in $C^J$. 
It has an inverse, say $r'$, in $C$ and since $1-r$ is invariant, so is $r'$. 
\end{pf}

Write $\bar{B}$ for $\hat{B}^I/\rad \hat{B}^I$. 
The lemma implies that if $\bar{B}^J$ is semisimple then $\hat{B}^H/\rad \hat{B}^H = \bar{B}^J$. 
We show that $\bar{B}^J$ is indeed semisimple. 
Note that the skew-commutation relation $b_\s b_\t = \z^{lp} b_\t b_\s$ shows that $b_\t$ restricts to isomorphisms $V_i \rtar V_{i+lp}$. 
Hence conjugation by $b_\t$ yields isomorphisms of various factors of $\bar{B}$, namely
\[ \t:k((s))^{d_{i+lp} \times d_{i+lp}} 
      \xrightarrow{\sim} k((s))^{d_i \times d_i} .\]
Now $\eta=\z^{lp}$ is a primitive $e$-th root of unity so $lp$ generates 
$p\Z$ modulo $ep$. Hence 
\[\left(\prod_{i=0}^{pe-1} k((s))^{d_i \times d_i}\right)^J 
  = \left(\prod_{i=0}^{p-1} k((s))^{d_i \times d_i}\right)^{\la\t^e\ra}  \]
We compute the invariants of each factor 
individually. They are symmetric so we assume $i =0$ 
and write $d= d_0$. Recall that $\t$ acts on $k((s))$ by 
$\t: s \mapsto \xi s$ where $\xi$ is a primitive $qe$-th root of unity.  
Also $b_{\t}^e$ has order $q$ so $b_\t^e|_{V_0}$ can be 
diagonalized to have block matrix form 
\begin{equation}  \label{bt} 
b_\t^e|_{V_0} = 
\begin{pmatrix}
I      &   0     &   0    &   \hdots  &  0 \\
0      & \xi^e I    &   0    &   \hdots  &  0 \\
0      &   0     & \ddots &           &  \vdots \\
\vdots & \vdots  &        &  \ddots   &  0      \\
0      &   0     & \hdots &    0      & \xi^{(q-1)e} I   
\end{pmatrix}
\end{equation}
It follows that 
\[(k((s))^{d \times d})^{\la\t^e\ra} = 
\begin{pmatrix}
k((s^q))   & s^{-1} k((s^q)) &  \hdots \\
sk((s^q))  &      k((s^q))   &  \hdots \\
\vdots     & \vdots          & \ddots
\end{pmatrix}
\simeq k((s^q))^{d \times d} .\]
Hence
\[ \hat{B}^H/\rad \hat{B}^H \simeq
 \prod_{i=0}^{p-1} k((s^q))^{d_i \times d_i} .\]
Normality is thus equivalent to 
all non-zero $d_i$ being equal. 
This gives the ramification in the statement. 
Finally, the index of the associated central simple algebra can be determined from the last equation. 
\end{pf}

The following lemma gives the condition for the order to be Gorenstein.  

\begin{lemma} \label{canonical}
Suppose that $G$ acts on $B = S \otimes k^{m \times m}$ as above. 
Let $\chi$ denote the character 
$G \subset GL_2 \xrightarrow{\det} k^*$. 
Suppose there exists an invertible matrix 
$\theta \in k^{m \times m}$ which lies in the $\chi^{-1}$ 
component of the $G$-module $k^{m \times m}$. Then 
$A= B^G$ is Gorenstein and the dualising sheaf can 
be identified as $A\theta = \theta A$. 
\end{lemma}
\begin{pf}
We first show that $\w_A=\w_B^G$. Let $R$ denote the centre of $A$.
Then $\w_A=\Hom_R(A,\w_R)$ and $\w_B=\Hom_S(B,\w_S)=\Hom_R(B,\w_R)$,
so $\w_B^G=\Hom_R(B,\w_R)^G=\Hom_R(B^G,\w_R)=\w_A$ as required. 

There is a natural isomorphism $B \otimes \omega_S  \rightarrow \omega_B$
which gives a $G$-module isomorphism $B \otimes_k V \rightarrow \omega_B$,
where $V=k \cdot du \wedge dv$, the $G$-module of dimension $1$ given by the 
character $\chi$. Hence $\w_A = \w_B^G$ is identified with the $\chi^{-1}$ 
component of $B$. As $\theta$ is invertible, 
left multiplication by $\theta$ defines a right $A$-module 
isomorphism $B^G \stackrel{\sim}{\rtar} \w_B^G$.
The argument is left-right symmetric, so $A$ is Gorenstein.
\end{pf}

\end{section}

\begin{section}{Quotient constructions}
\label{exp}
In this section we give explicit quotient constructions for canonical singularities of orders. 
The analysis is case-by-case, depending on the ramification data. 
We deduce that canonical orders are Gorenstein and satisfy a basic version of the McKay correspondence. 

We use the notation of sections~\ref{frt} and \ref{sinvariant}. 
Let $A$ be a canonical order with centre $R$ a complete local normal domain of dimension $2$, 
whose ramification data is of some fixed type as classified in Thm.~\ref{canramdata}. 
We first find a subgroup $G \subset GL(ku+kv)$ acting on $S=k[[u,v]]$ such that $S/R$ has the same ramification as the order $A$. 
Then Thm.~\ref{artincover} shows that $A = B^G$ for some lift of the $G$-action to $B = S^{m \times m}$. 
As noted in Lem.~\ref{projectiverep}, the $G$-action on $B$ is determined by a projective representation $\b \colon G \rightarrow \PGL_m$ 
which in turn is given by a representation $b$ of some central extension $G'$ of $G$ by a finite cyclic group $\mu_d$. 

Using Thm.~\ref{normal}, we compute which representations $b$ of $G'$ yield normal orders with the correct ramification data. 
We then verify that there is an invertible matrix satisfying the condition of Lem.~\ref{canonical} and deduce that $A = B^G$ is Gorenstein. 

In Prop.~\ref{arquivers}, we showed that reflexive $A$-modules correspond to $G'$-modules where the central subgroup $\mu_d$ acts via some 
fixed character. 
Such a representation $b$ of $G'$ is said to be {\it permissible} if it satisfies the following: 
\begin{enumerate}
\item For every line of reflection $L$ of $G$, let $\s' \in G'$ be a lift of a generator of the inertial group of $L$. 
Then the eigenspaces of $b_{\s'}$ all have the same dimension (this is the eigenspace condition in Thm.~\ref{normal}).
\item $b$ is irreducible amongst such representations.
\end{enumerate}
We also call the corresponding $A$-modules {\it permissible}. 
We verify below that the number of such permissible modules is one more than the number of exceptional curves in a minimal resolution. 

Let $K$ be the field of fractions of $R$ and $A_K:=A \otimes_R K$ the central simple algebra associated to $A$. 

\begin{subsection}{Type $A_{1,2,\xi}$}   \label{typebe}

Let $G \simeq \Z/2e \times \Z/2e$ and $\s,\t$ be the generators of the two cyclic groups. 
Let $\z$ be a primitive $2e$-th root of unity and suppose the group acts on $S$ by 
$\s:u \mapsto \z u, v \mapsto v$ and $\t:u \mapsto u, v \mapsto \z v$. 
Then $R:= S^G = k[[x,y]]$ where $x = u^{2e}, y = v^{2e}$ and $S/R$ is a ramified cover which is ramified along $xy = 0$ with ramification index $2e$.
This coincides with the ramification of a type $A_{1,2,\xi}$ canonical order $A$. 
Note that the stabiliser of $(u=0)$ and $(v=0)$ is $G$ and the inertia groups of $(u=0)$ and $(v=0)$ are $\la\t\ra$ and $\la\s\ra$ respectively. 

We compute the possibilities for the projective representation $\beta:G \rtar \PGL_m$. 
As usual, we lift $\beta$ to a representation $b$ of $G'$. 
We let $\beta_\s$ be conjugation by $b_\s \in \GL_m$ and $\beta_\t$ be conjugation by $b_\t \in \GL_m$. 
Since $\b_\s^{2e} = 1, \b_\t^{2e} = 1$, we may scale so that $b_\s^{2e} = 1, b_\t^{2e} = 1$. 
(To do this we may need to enlarge $G'$ as described in Rmk.~\ref{centralextensionremark}.) 

The cyclic covers of the components of the discriminant have ramification index $e$ over the closed point. 
So, by Thm.~\ref{normal}, we have $\xi b_\t b_\s = b_\s b_\t$ for some primitive $e$-th root of unity $\xi$. 
This is the $\xi$ in the subscript of $A_{12,\xi}$. 
Write $\xi = \z^{2l}$ where $l$ is relatively prime to $e$. Thus
\begin{equation}  \label{reln}
b_\s^{2e} =  b_\t^{2e} = 1 \ , \ 
    \z^{2l} b_\t b_\s = b_\s b_\t.
\end{equation}
The remaining condition required by Thm.~\ref{normal} is that both $b_{\s},b_{\t}$ have all $2e$ possible eigenvalues and that 
all the eigenspaces have the same dimension. 

We let $V$ be the $G'$-module corresponding to $b$.
We first describe all irreducible $G'$-modules which satisfy (\ref{reln}) and 
second determine what combinations of these satisfy the eigenspace condition of Thm.~\ref{normal}. 

We introduce the following $e \times e$-matrices. 
\[ P = 
\begin{pmatrix}
1      &   0     & \hdots &  0 \\
0      & \z^2    & \ddots &   \\
\vdots & \ddots  & \ddots &  0      \\
0      &         & 0      & \z^{2(e-1)}   
\end{pmatrix}  \ 
Q_{\pm} = 
\begin{pmatrix}
0      & \hdots   & 0      &  \pm 1  \\
1      &      0   &        &  0      \\
0      &      1   & \ddots &         \\
\vdots & \ddots   & \ddots & \ddots  & 
\end{pmatrix}  \ 
N = 
\begin{pmatrix}
1      &   0     & \hdots &  0 \\
0      & \z^{-1}    & \ddots &   \\
\vdots & \ddots  & \ddots &  0      \\
0      &         & 0      & \z^{1-e}   
\end{pmatrix}  \]

Let $W$ be an irreducible $G'$-module. 
Suppose the corresponding representation, which we also denote by $b$, satisfies the equations (\ref{reln}). 
We consider the eigenspace decomposition $W = \oplus_{i \in \Z/2e} W_i$ with respect to $b_\s$ where $W_i$ is the $\z^i$-eigenspace of $b_\s$. 
From (\ref{reln}) we see that $b_\t$ restricts to isomorphisms $W_i \rtar W_{i+2l}$. 
Now $b_\t^e$ maps $W_i \rtar W_i$ and $(b^e_{\t})^2=1$. 
Suppose that $W_i$ is non-zero and pick a $b_\t^e$-eigenvector $w \in W_i$ which must have eigenvalue $\pm 1$. 
We may assume that $i = 0$ or 1. 
Irreducibility shows that $W$ has basis $\{w, b_\t w, \ldots , b_\t^{e-1} w\}$. 
We write $W_{i\pm}$ for this $G'$-module where $b_\t = Q_{\pm}$ and $b_\s = P^l$ if $i = 0$ and $b_\s = \z P^l$ if $i = 1$. 

The eigenspace condition of Thm.~\ref{normal} shows that $A = B^G$ is a canonical order of type $A_{12,\xi}$ precisely when $V$ is a direct sum of 
modules of the form $W_{0+} \oplus W_{1-}$ and $W_{0-} \oplus W_{1+}$. 
To prove that $A$ is Gorenstein, we may assume that $V$ is one of these two modules. 
By Lem.~\ref{canonical}, $A$ is Gorenstein if there is an invertible matrix $\theta \in k^{m \times m}$ such that 
$b_\s^{-1}\theta b_\s = \z^{-1} \theta,   b_\t^{-1}\theta b_\t = \z^{-1} \theta$.
In the first case, we have 
$ b_\s = 
\bigl(\begin{smallmatrix}
 P^l & 0 \\ 0 & \z P^l
\end{smallmatrix}\bigr),
b_\t = 
\bigl(\begin{smallmatrix}
 Q_+ & 0 \\ 0 & Q_-
\end{smallmatrix}\bigr)$. 
Let $j$ be the inverse of $l$ modulo $e$, then set
$\theta = 
\bigl(\begin{smallmatrix}
 0 &  Q_+^j N \\ N & 0 
\end{smallmatrix}\bigr)$. 
In the second case we have, 
$ b_\s = 
\bigl(\begin{smallmatrix}
 P^l & 0 \\ 0 & \z P^l
\end{smallmatrix}\bigr),
b_\t = 
\bigl(\begin{smallmatrix}
 Q_- & 0 \\ 0 & Q_+
\end{smallmatrix}\bigr)$
and we set
$\theta = 
\bigl(\begin{smallmatrix}
 0 & Q_-^j N \\ N & 0
\end{smallmatrix}\bigr)$. 

Note that there are exactly two permissible representations ($W_{0+} \oplus W_{1-}$ and $W_{0-} \oplus W_{1+}$ ) 
which is one more than the number of exceptional curves in a minimal resolution. 

\end{subsection}

\begin{subsection}{Type $BL_n$}\label{BLn}

Let $r = 2n+1$ and let $G$ be the dihedral group 
$\la \s,\t \st \s^r = \t^2 = 1, \s\t=\t\s^{-1} \ra$. 
Suppose $G$ acts linearly on $S = k[[u,v]]$ by 
\[ \s = 
\begin{pmatrix}
 \z & 0 \\
  0 & \z^{-1} 
\end{pmatrix} \ \ , \ \ 
  \t = 
\begin{pmatrix}
  0 & 1 \\
  1 & 0
\end{pmatrix} 
\]
where $\z$ is a primitive $r$-th root of unity. 
We compute $R := S^G = k[[x,y]]$ where $x = uv, y = \frac{1}{2}(u^r + v^r)$. 
Ramification of $S/R$ occurs at the fixed lines of the pseudo-reflections of $G$. 
There is one conjugacy class of pseudo-reflections so the ramification curve is the image in $\Spec R$ of $(\Spec S)^{\la\psi\ra}$ 
where $\psi$ is any pseudo-reflection. 
Picking $\psi = \t$, we see that $S/R$ is ramified over the image of $(u=v)$ which is the cuspidal curve $(y^2 = x^r)$. 
The ramification index is $2$ so this is the same as the ramification of a canonical order of type $BL_n$. 

For this type, $A_K$ is trivial in the Brauer group, so the projective representation $\b$ 
lifts to an actual representation $b$. 
We consider the irreducible one-dimensional representations of $G$ 
\[ \r^0:\s,\t \mapsto 1 \ \ , 
    \ \ \r^-:\s \mapsto 1, \t \mapsto -1  \]
and the irreducible two-dimensional representations 
\begin{equation}  \label{twodim} 
\r^i:\s \mapsto 
\begin{pmatrix}
 \z^i  &  0 \\
    0  & \z^{-i}
\end{pmatrix},
\t \mapsto 
\begin{pmatrix}
   0 & 1 \\
   1 & 0
\end{pmatrix}
\end{equation}
where $i= 1, \ldots ,n$. 

The inertia group at $(u=v)$ is $\la\t\ra$ so the eigenspace condition of Thm.~\ref{normal} shows that $b$ is 
the direct sum of $\rho^i$ and $\rho^0 \oplus \rho^-$. 
Note that again we have $n+1$ permissible representations, namely the $\rho^i$ and $\rho^0 \oplus \rho^-$.  
To show $A$ is Gorenstein in this case, we may assume that $b$ is one of these. 

The determinant character $\chi \colon G \rtar k^*$ is given by $\s \mapsto 1$, $\t \mapsto -1$. 
If $b = \rho^i$ then setting 
$\theta = 
\bigl(
\begin{smallmatrix}
1 & \\  & -1 
\end{smallmatrix}
\bigr)$ 
in Lem.~\ref{canonical} shows that $A$ is Gorenstein. 
If $b = \rho^0 \oplus \rho^-$  then we set 
$\theta = 
\bigl(
\begin{smallmatrix}
 & 1 \\ 1 &  
\end{smallmatrix}
\bigr)$.

\end{subsection}

\begin{subsection}{Type $B_n$}\label{Bn}

Let $r = 2n$ and let $G$ be the dihedral group 
$\la \s,\t \st \s^r = \t^2 = 1, \s\t = \t\s^{-1} \ra$. 
Suppose $G$ acts linearly on $S = k[[u,v]]$ as in Sec.~\ref{BLn} (except now $r$ is even). 
As before, we find $R:= S^G = k[[x,y]]$ where $x = uv, y = \frac{1}{2}(u^r + v^r)$. 
The difference is that there are two conjugacy classes of pseudo-reflections, namely, 
$\{ \s^i\t \st i \ \mbox{even}\}$ and $\{ \s^i\t \st i \ \mbox{odd}\}$. 
Let $L=(u = v)$ be the line fixed by $\t$ and $L'=(u = \z^{-1} v)$ the line fixed by $\s\t$. 
Their images in $\Spec R$ give the irreducible components of the ramification locus, namely, $(y = x^n)$ and $(y = -x^n)$. 
Furthermore, the pointwise stabilisers of $L,L'$ are $\la\t\ra,\la\s\t\ra$. 
These have order $2$ so $S/R$ ramifies with ramification index $2$. 
This coincides with the ramification of a canonical order $A$ of type $B_n$. 

In this case, as before, $A_K$ is trivial in the Brauer group so $\b \in H^1(G,\PGL_m)$ lifts to an actual representation $b \in H^1(G,\GL_m)$. 
Again we decompose $b$ into irreducible representations. 

As in Sec.~\ref{BLn}, there are irreducible two-dimensional representations $\r^i$, for $i = 1, \ldots, n-1$ defined by (\ref{twodim}). 
There are however, $4$ irreducible one-dimensional representations.
\begin{align*}
\r^{00}: \s \mapsto 1, \t \mapsto 1  
\hspace{1cm}&
\r^{01}: \s \mapsto 1, \t \mapsto -1  \\
\r^{10}: \s \mapsto -1, \t \mapsto 1  
\hspace{1cm}&
\r^{11}: \s \mapsto -1, \t \mapsto -1  
\end{align*}
We compute $\Stab L =\la\t,\s^n\ra, \Stab L' = \la\s\t,\s^n\ra$ 
and the inertia groups at $L,L'$ are $\la\t\ra, \la\s\t\ra$ respectively. 
By taking direct sums of matrices and using the eigenspace condition of Thm.~\ref{normal}, we deduce 
$b = \r^i, \r^{00} \oplus \r^{01}$ or $\r^{10} \oplus \r^{11}$. 
Hence again we see that there are $n+1$ permissible 
representations.

The determinant character $\chi:G \rtar k^*$ is given by $\s \mapsto 1$, $\t \mapsto -1$. 
For $b = \r^i$, setting 
$\theta =  \bigl(
\begin{smallmatrix}
 1 & 0 \\
 0 & -1
\end{smallmatrix} 
\bigr)$ 
shows that $A$ is Gorenstein. 
In the other two cases, we use bases for $\r^{00} \oplus \r^{01}$ and $\r^{10} \oplus \r^{11}$ which are compatible with the direct sum decomposition. 
Then setting 
$\theta =  \bigl(
\begin{smallmatrix}
 0 & 1 \\
 1 & 0
\end{smallmatrix} 
\bigr)$ 
shows that $A$ is Gorenstein. 

\end{subsection}

\begin{subsection}{Type $L_n$}\label{Ln}

This is Artin's type $II_k$ of \cite{A86}.

The ramification curves and the ramification indices for types $B_{n+1}$ and $L_n$ coincide so if $G$ is the dihedral group of Sec.~\ref{Bn} 
with $r=2n+2$, and $G$ acts on $S$ as before, then $S/R$ has the same ramification as a type $L_n$ canonical order $A$.

The difference in this case is that $A_K$ is non-trivial in the Brauer group. 
Write $\b$ for the projective representation describing the group action on $B$ and lift $\b_\s,\b_\t$ to elements $b_\s,b_\t \in \GL_m$. 
We may assume  
\begin{equation} \label{relf1n}
b_\s^{2n+2} = b_\t^2 = 1 \ , \ 
 b_\t b_\s^{-1} = \l b_\s b_\t  
\end{equation}
for some scalar $\l$. 

The cyclic covers of the components of the discriminant ramify over the closed point.
So, using the stabiliser groups computed in Sec.~\ref{Bn} and Thm.~\ref{normal}, we deduce the relation 
$b_\s^{n+1} b_\t = -b_\t b_\s^{n+1}$, 
or, equivalently, $\l^{n+1} = -1$. 
Write $\l = \z^{-a}$ where $\z$ is a primitive $(2n+2)$-th root of unity, and note that $a$ is odd. 
(Note also that $b^{n+1}_\s$ switches the $(\pm 1)$-eigenspaces of $b_\t$ so normality imposes no conditions as long as $\l^{n+1} = -1$.)

We view $b$ as a representation of $G'$ as usual and let $V$ be the corresponding $G'$-module. 
To prove that $A$ is Gorenstein using Lem.~\ref{canonical}, we may assume that $b$ is irreducible and still satisfies (\ref{relf1n}). 
We first decompose $V$ into $\z^i$-eigenspaces $V_i$ with respect to $b_\s$. 
Note that there are induced isomorphisms $b_\t:V_i \xrightarrow{\sim} V_{a-i}$. 
Furthermore, $a$ being odd implies that these eigenspaces are distinct. 
Using (\ref{relf1n}), we see that the irreducible representations have the form 
\[ b_\s = 
\begin{pmatrix}
 \z^i & 0 \\
   0  & \z^{a-i} 
\end{pmatrix} \ \ , 
\ \ b_\t = 
\begin{pmatrix}
   0  & 1 \\
   1  & 0 
\end{pmatrix} .\]
Setting 
$\theta = 
\bigl(
\begin{smallmatrix}
1 & 0 \\
0 & -1
\end{smallmatrix}
\bigr)$ in Lem.~\ref{canonical} shows that 
$A$ is Gorenstein. Finally, since swapping 
eigenspaces for $\z^i$ and $\z^{a-i}$ gives 
isomorphic modules, there are $r/2 = n+1$ 
permissible modules. 

\end{subsection}

\begin{subsection}{Type $DL_n$}
This is Artin's type $III_k$ of \cite{A86}. 

Let $G$ be the subgroup of $\GL_2$ generated by 
\[ \s = 
\begin{pmatrix}
 \z & 0 \\
  0 & \z^{-1}
\end{pmatrix} \ \ , 
\ \ \t = 
\begin{pmatrix}
 0 & 1 \\
  1 & 0
\end{pmatrix} \ \ , 
\ \ \pi = 
\begin{pmatrix}
 1 & 0 \\
  0 & -1
\end{pmatrix} \]
where $\z$ is a primitive $2r$-th root of unity and $r=2n-1$ is an odd integer. 
Consider the linear action of $G$ on $S$ defined by the above matrices, and let $R=S^G$. 
We compute that $S/R$ has the same ramification as a canonical order $A$ of type $DL_n$. 
Note first that $R = k[[x,y]]$ where $x = u^2v^2, y = \frac{1}{2}(u^{2r}+ v^{2r})$. 
There are two conjugacy classes of pseudo-reflections, $\{\pi,\s^r\pi\}$ and $\{\s^i\t\}$. 
The fixed line $L_0=(u = 0)$ of $\pi$ corresponds to the discriminant curve $(x=0)$, 
the fixed line $L_1=(u = v)$ of $\t$ corresponds to the discriminant curve $(y^2 = x^r)$, 
and in each case the ramification index is $2$.

As usual, we consider the projective representation $\b$ defining $A$. 
We lift $\b_\s,\b_\t,\b_\pi$ to $b_\s,b_\t,b_\pi \in \GL_m$ such that $b_\s^{2r} = b_\t^2 = b_\pi^2 = 1$. 
The cyclic covers of the components of the discriminant ramify over the closed point, and $\Stab L_0 = \la\s,\pi\ra, \Stab L_1 = \la\t, \s^r\ra$, 
so Thm.~\ref{normal} implies  
\begin{equation} \label{cond5}
b_\s b_\pi = - b_\pi b_\s \ , \ 
   b_\s^r b_\t = - b_\t b_\s^r  .  
\end{equation}
As in Sec.~\ref{Ln}, any $b$ satisfying the above equations gives a canonical order of type $DL_n$. 

To show that $A$ is Gorenstein using Lem.~\ref{canonical}, we may assume $b$ corresponds to an irreducible representation $V$ of $G'$. 
Note first that $b_\s b_\t = \l b_\t b_\s^{-1}$ for some scalar $\l$ and (\ref{cond5}) gives $\l^r = -1$. 
Hence $\l = \z^a$ for some odd integer $a$. 
Also, $b_{\pi}b_{\t} = \rho b_\s^r b_\t b_\pi$ for some scalar $\rho$. 
(It turns out that $\rho^2 = -1$ though we do not need this.) 
We decompose $V$ into $\z^i$-eigenspaces $V_i$ with respect to $b_\s$. 
The relations (\ref{cond5}) show that there are induced isomorphisms $b_\t:V_i \rtar V_{a-i}$, $b_\pi:V_i \rtar V_{r+i}$. 
Suppose that $w$ is a non-zero vector in $V_i$ and that $b_\t w,b_\pi w$ are linearly independent. 
Then, since $b_\t^2=b_\pi^2=1$, $V$ has basis $\{ w,b_\pi w, b_\t b_\pi w, b_\t w\}$ and 
\[ b_\s = 
\begin{pmatrix}
 \z^i &  &  &  \\
      & -\z^i & & \\
     &   & -\z^{a-i} & \\
  &  &  &  \z^{a-i}
\end{pmatrix} \ , \ 
b_\t = 
\begin{pmatrix}
  &  &  & 1 \\
      &  & 1  & \\
     & 1  &  & \\
  1 &  &  & 
\end{pmatrix} \ , \ 
b_\pi = 
\begin{pmatrix}
  & 1 &  &  \\
    1  &  & & \\
     &   &  & \pm \rho \\
  &  &  \mp \rho & 
\end{pmatrix}
\]
where the sign in $b_\pi$ depends on $(-1)^i$. 
Setting 
$$\theta = 
\left(
\begin{smallmatrix}
1  &  &  &  \\
  & -1 &  &  \\
  &  & 1 &  \\
  &  &  & -1 
\end{smallmatrix}
\right)$$ 
in Lem.~\ref{canonical} shows that $A$ is Gorenstein. 
Furthermore, the module is an irreducible $G'$-module unless $i$ satisfies $a-i \equiv r+i \mod 2r$, i.e., $2i \equiv a-r$, 
in which case the module decomposes into two 2-dimensional modules. 
In the irreducible case, the eigenvalues of $b_\s$ play a symmetric role so swapping $i$ with $a-i, i+r$ or $a-i+r$ gives an isomorphic module. 
There are consequently $n-1$ permissible modules of dimension $4$. 

For the two dimensional modules, we have $b_\pi w = \nu b_\t w$ for some scalar $\nu$, so that $V$ has basis $\{w,b_\t w\}$ and  
\[ b_\s = 
\begin{pmatrix}
 \z^i &    \\
      & -\z^{i} 
\end{pmatrix} \ , \ 
b_\t = 
\begin{pmatrix}
     & 1   \\
  1 & 
\end{pmatrix} \ , \ 
b_\pi = 
\begin{pmatrix}
  & (-1)^i \rho \nu   \\
    \nu  & 
\end{pmatrix}.
\]
In this case, setting 
$\theta = 
\bigl(
\begin{smallmatrix}
1  &    \\
  & -1 
\end{smallmatrix}
\bigr)$ 
in Lem.~\ref{canonical} shows that $A$ is Gorenstein. 
Note that as $b_{\pi}^2 = 1$, there are exactly two possible choices for $\nu$. 
If $i$ satisfies $2i \equiv a-r$ then the other value which satisfies the congruence is $i+r$, and changing $i$ to $i+r$ yields no new modules. 
Hence there are two permissible modules of dimension $2$ giving a total of $n+1$ permissible modules. 

\end{subsection}

\begin{subsection}{Type $BD_n$}

For convenience, we set $p=n-1$, and let $\z$ be a primitive $4p$-th root of unity. 
Let $G$ be the subgroup of $\GL_2$ generated by 
\[ \s = 
\begin{pmatrix}
 \z & 0 \\
  0 & \z^{-1}
\end{pmatrix} \ \ , 
\ \ \t = 
\begin{pmatrix}
  0 & 1 \\
 -1 & 0
\end{pmatrix} \ \ , 
\ \ \pi = 
\begin{pmatrix}
 -1 & 0 \\
  0 & 1
\end{pmatrix} \]
Note the following relations 
\[ \s^{4p} = \pi^2 = 1, \t^2 = \s^{2p}, 
  \s \t = \t \s^{-1}, \t \pi = \pi \t^{-1}, 
  \pi \s = \s \pi                    .  \]
Let $G$ act linearly on $S= k[[u,v]]$ via the matrices above and set $R = S^G$. 
We claim that $S/R$ has the same ramification as a canonical order of type $BD_n$. 
Note first that $R = k[[x,y]]$ where $x = u^2v^2, y =\frac{1}{2}(u^{4p}+v^{4p})$. 
There are 3 conjugacy classes of pseudo-reflections $\{\pi,\s^{2p}\pi\}$, $\{\s^i \t\pi \st i\ \mbox{odd} \}$, and $\{\s^i \t\pi \st i\ \mbox{even} \}$. The lines fixed by the three representative pseudo-reflections $\pi, \t\pi, \s\t\pi$ are 
$$ L_0=(u = 0) \ , \ L_1=(u = v)  \ , \ L_2=(u = \z v)      .$$
The images of these lines in $\Spec R$ are the discriminant curves  
$$ (x = 0) \ , \ (y = x^{p}) \ , \ (y = -x^{p})$$
of $S/R$.
The ramification index equals $2$ in each case so the ramification is the same as that of a canonical order of type $BD_n$. 
Moreover, for such an order $A$, we may assume that the cyclic covers of $(x=0)$ and $(y=x^p)$ determined by $A$ ramify over the closed point, 
while the cyclic cover of $(y=-x^{p})$ is unramified.

Again, let $\b$ be the projective representation of $G$ which gives the $G$-action on $B$ and 
$b$ its lift to an actual representation of $G'$. 
We write $b_\s,b_\t,b_\pi$ as before. 
To simplify calculations, we normalize $b_\s, b_\t, b_\pi$ so that $b_\s^{4p} = b_\pi^2 = 1 , b_\t^2 = b_\s^{2p}$. 

Now $\Stab L_0 = \la\s,\pi\ra$ so Thm.~\ref{normal} shows that $b_\s b_\pi = - b_\pi b_\s$. 
Suppose now that $b_\s b_\t = \l b_\t b_\s^{-1} , b_\t b_\pi = \mu b_\pi b_\t^{-1}$. 
Since $b_\pi$ commutes with $b_\t^2 = b_\s^{2p}$, we must have $\mu^2 = 1$. 
Also, commuting $b_\s^{2p} = b_\t^2$ through $b_\t$ shows that $\l^{2p} = 1$. 

We have $\Stab L_1 = \la\t \pi, \s^p \pi\ra$, so Thm.~\ref{normal} yields  $b_\t b_\pi b_\s^p b_\pi = - b_\s^p b_\pi b_\t b_\pi$, 
which amounts to $\mu \l^p = (-1)^{p+1}$. 
Also $\Stab L_2 = \la\s\t\pi , \s^p \pi\ra$, which gives the commutativity relation $b_\s b_\t b_\pi b_\s^p b_\pi =b_\s^p b_\pi b_\s b_\t b_\pi$ 
because the cyclic cover of the image of $L_2$ is unramified.
This yields the same condition $\mu \l^p = (-1)^{p+1}$ as before. 
We also have the eigenspace condition on $b_\s b_\t b_\pi$. 
Since $(b_\s b_\t b_\pi)^2$ is a scalar, this amounts to the trace of $b_\s b_\t b_\pi$ being zero. 

We assume as usual that $b$ is irreducible and let $V$ be the corresponding $G'$-module. 
We begin by computing all irreducible $G'$-modules where $\mu_d = \ker (G' \rtar G)$ acts appropriately. 
Let $V = \oplus V_i$ be the eigenspace decomposition with respect to $b_\s$ where $V_i$ has eigenvalue $\z^i$. 
Write $\l = \z^a$ and observe that $a$ is even since $\l^{2p} =1$.  
Note that we have induced isomorphisms $b_\t: V_i \rtar V_{a-i}, b_\pi:V_i \rtar V_{2p+i}$. 
Suppose first that we have a non-zero element $w \in V_i$ and that $w,b_\pi w, b_\t b_\pi w, b_\t w$ are linearly independent. 
Taking into account how $b_\t,b_\pi$ permute the $V_i$ we find 
\[ b_\s = 
\begin{pmatrix}
 \z^i &  &  &  \\
      & -\z^i & & \\
     &   & -\z^{a-i} & \\
  &  &  &  \z^{a-i}
\end{pmatrix} \]
\[ b_\t = 
\begin{pmatrix}
  &  &  & (-1)^i \\
      &  & (-1)^i  & \\
     & 1  &  & \\
  1 &  &  & 
\end{pmatrix} \ , \ 
b_\pi = 
\begin{pmatrix}
  & 1 &  &  \\
    1  &  & & \\
     &   &  & \mu (-1)^i \\
  &  &  \mu (-1)^i & 
\end{pmatrix}.
\]
One checks easily that the eigenspace condition of Thm.~\ref{normal} is satisfied and that setting 
\[ \theta = 
\begin{pmatrix}
1 & & & \\
 & -1 & & \\
 & & -1 & \\
 & & & 1
\end{pmatrix}\]
in Lem.~\ref{canonical} shows that $A$ is Gorenstein. 

This $G'$-module is irreducible except when either (1) $i \equiv a-i \mod 4p$, or (2) $a-i \equiv 2p+i \mod 4p$. 
In these cases, the module decomposes into a direct sum of two 2-dimensional modules. 
As in the type $DL_n$ case, we see there are $p-1$ permissible modules of dimension $4$. 

Assume now that we are in case (1) and that say $w,b_\t w$ are linearly dependent. 
Note that $b_\t^2 = b_\s^{2p}$ forces $b_\t w = \pm \z^{ip} w$ in this case. 
Using the basis $w,b_\pi w$ we find 
\[ b_\s = 
\begin{pmatrix}
 \z^i &    \\
      & -\z^i 
\end{pmatrix} \ , \ 
b_\t = 
\begin{pmatrix}
  \pm \z^{ip}  &    \\
    & \pm \mu\z^{-ip}
\end{pmatrix} \ , \ 
b_\pi = 
\begin{pmatrix}
  &  1   \\
    1  & 
\end{pmatrix}.
\]
The trace of $b_\s b_\t b_\pi$ is zero so the eigenspace condition of Thm.~\ref{normal} is satisfied. 
Setting 
$\theta = \bigl(
\begin{smallmatrix}
 1 & \\ & -1
\end{smallmatrix}\bigr)$ 
in Lem.~\ref{canonical} shows that $A$ is Gorenstein. 
Note that there is a choice of sign in $b_\t$ and that swapping values of $i$ corresponds to switching the roles of $w,b_\pi w$. 
Hence there are two new permissible modules. 

Finally, suppose we are in case (2). 
We let $b_\t w = \nu b_\pi w$ for some scalar $\nu$. 
The module $V$ will sometimes be denoted $V^\nu$ to emphasize the dependence on $\nu$. 
Note first that if $i$ is a solution to $a-i \equiv 2p+i \mod 4p$ then the only other solution is $i + 2p$. 
Since $V = V_i \oplus V_{i+2p}$, the latter does not give any new $G'$-modules. 
The fact that $a-2p \equiv 2i \mod 4p$ and our ramification condition on $\mu, \l$ give 
\[ \mu = \l^{-p}(-1)^{p+1} = \z^{-ap + 2p^2 + 2p} = 
   \z^{-2ip + 2p} = -\z^{-2ip}.\]
Hence, 
\[ b_\t b_\pi w = \mu b_\pi b_\t^{-1} w = 
   \mu b_\pi b_\t b_\s^{2p} w = \mu \z^{2ip} \nu w 
                                      = -\nu w .\]
We can now compute $V$ as 
\[ b_\s = 
\begin{pmatrix}
 \z^i &    \\
      & -\z^i 
\end{pmatrix} \ , \ 
b_\t = 
\begin{pmatrix}
   & -\nu   \\
  \nu  & 
\end{pmatrix} \ , \ 
b_\pi = 
\begin{pmatrix}
  &  1   \\
    1  & 
\end{pmatrix}. \]
There are precisely two possible values for $\nu$ since 
$b_\s^{2p} = b_\t^2$. 
Now the trace of $b_\s b_\t b_\pi$ is non-zero and the eigenspace condition of Thm.~\ref{normal} fails. 
To obtain the eigenspace condition, we see that the multiplicity of $V^{\nu}$ in $V$ must be the same as the multiplicity of $V^{-\nu}$ in $V$. 
We may thus suppose that 
$V = V^{\nu} \oplus V^{-\nu}$ so that 
\[ b_\s = 
\begin{pmatrix}
 \z^i &  &  &  \\
      & -\z^i & & \\
     &   & \z^i & \\
  &  &  &  -\z^i
\end{pmatrix} \ , 
\ b_\t = 
\begin{pmatrix}
  & -\nu &  &  \\
 \nu  &  &   & \\
     &   &  & \nu \\
   &  & -\nu & 
\end{pmatrix} \ , \ 
b_\pi = 
\begin{pmatrix}
  & 1 &  &  \\
    1  &  & & \\
     &   &  & 1 \\
  &  &  1 & 
\end{pmatrix}.
\]
Setting 
\[ \theta = 
\begin{pmatrix}
 & & 1 & \\
 & & & -1 \\
1 & & & \\
 & -1 & & 
\end{pmatrix} .\]
in Lem.~\ref{canonical} shows that $A$ is Gorenstein in this case too. 
This gives one more permissible module for a total of $p+2 = n+1$. 
\end{subsection}

\begin{subsection}{Type $ADE$}
Let $A$ be an order of type $A$, $D$ or $E$ so that its centre is of the form $R = S^G = k[[u,v]]^G$ where $G$ is a finite subgroup of $SL_2$. 
Let $B\simeq S^{m \times m}$ be the Artin cover of $A$ with respect to $S/R$.
Then $A = B^G$ is Gorenstein by Lem.~\ref{canonical} because $G$ acts linearly on $S$ by matrices of determinant $1$.

\end{subsection}

\begin{subsection}{Type $A_{n,\xi}$} 
Set $p = n+1$. Let $G$ be the subgroup of $GL_2$ generated by 
\[ \s = 
\begin{pmatrix}
 \z & 0 \\
  0 & \z^{-1} 
\end{pmatrix} \ \ , \ \ 
  \t = 
\begin{pmatrix}
  1 & 0 \\
  0 & \z^p
\end{pmatrix} 
\]
where $\z$ is a primitive $pe$-th root of unity. 
Note that $R:=S^G = k[[u^{pe},u^ev^e,v^{pe}]]$ is the rational double point of type $A_n$. 
In fact, if we let $H \leq G$ be the subgroup generated by $\s^p,\t$, then $S^H=k[[u^e,v^e]]$ is the smooth cyclic cover of $R$ unramified 
away from the singularity and $S/S^H$ is the $\Z/e \times \Z/e$-cover ramified above the normal crossing lines $(u^ev^e = 0)$. 
Consequently, $S/R$ has the same ramification as a canonical order $A$ of type $A_{n,\xi}$.

On $\Spec S$, the cover $S/R$ ramifies on the two lines 
$(u=0)$, $(v=0)$. Consider the line $L=(v=0)$, which is stabilised by the whole group $G$. 
The inertia group in this case is $I = \langle \t \rangle$.
As usual, let $b:G' \rtar GL_m$ be the representation of $G'$ corresponding to the canonical order $A$. 
Since the cyclic covers of the components of the discriminant ramify with index $e$ over the closed point, 
we have $b_{\s}b_{\t} = \xi b_{\t}b_{\s}$ for some primitive $e$-th root of unity. 
This is the $\xi$ in the subscript of $A_{n,\xi}$. 
We can write $\xi = \z^{lp}$ where $l$ is relatively prime to $e$. 
Let $V$ be the $G'$-module corresponding to $b$, and assume $V$ is irreducible. 
Let $V = \oplus V_i$ be the decomposition of $V$ into $\z^i$-eigenspaces with respect to $b_{\s}$. 
Now $b_{\t}$ induces isomorphisms $V_i \rtar V_{i+lp}$ and $b_{\t}^e = 1$ so, with respect to an appropriate basis of eigenvectors, we have 
$b_{\s} = \z^iN^{-lp}, b_{\t} = Q_+$ using the notation of Sec.~\ref{typebe}. 
Picking $j$ so that $jl \equiv 1$ modulo $e$ and setting $\theta = N^{jp}$, we see that $A$ is Gorenstein. 
Up to rearranging the eigenspaces, there are $p = n+1$ choices for the value of $\z^i$ and hence $n+1$ permissible modules. 

\end{subsection}

\begin{subsection}{A non-Gorenstein log terminal order} \label{ltgor}
Consider normal orders with centre $k[[x,y]]$, ramified on normal crossing lines with ramification index 3 
but where the cyclic covers are unramified at the node. 
We describe two examples of such orders, one of which is Gorenstein while the other is not. 
Consequently the Gorenstein property is not determined by the ramification data. 

Let $G$ be the group $\Z/3 \times \Z/3$ and $\s,\t$ the generators of the two cyclic groups. 
Let $\z$ be a primitive cube root of unity and suppose the group $G$ acts on $S=k[[u,v]]$ by 
$\s:u \mapsto \z u, v \mapsto v$ and $\t:u \mapsto u, v \mapsto \z v$ so that $R:= S^G = k[[x,y]]$ where $x = u^3, y = v^3$. 
Note $S/R$ is ramified along $xy = 0$ with ramification index 3.

We let $G$ act on $B:=S^{3 \times 3}$ via the representation $b: G \rtar GL_3$ defined by 

\[ b_{\s} = \begin{pmatrix} 1 & & \\ & \z & \\ & & \z^2 
                                           \end{pmatrix}, 
   b_{\t} = \begin{pmatrix} 1 & & \\ & \z^2 & \\ & & \z 
                                           \end{pmatrix} .\]
Let $A$ be the order $B^G$ which we note is normal by Thm.~\ref{normal}. 
Note also that $A$ is log terminal by Thm.~\ref{tltisfrt}. 

We view $S$ as the completion of a graded ring which is graded by degree in $u,v$. 
Note that the action of $G$ on $B$ is graded and, identifying $\omega_B$ with $B$, we see that its action on $\w_B$ is also graded. 
The degree zero component of $A:= B^G$ is $A_0 = k^3$, the set of diagonal matrices over $k$. 
The lowest degree component of $\w_B^G$ is degree 1 and it is free of rank two over $A_0$. 
Indeed it is generated by 
\[ \begin{pmatrix}  & v & \\ &  & v \\ v& & 
                                           \end{pmatrix}, 
   \begin{pmatrix}  & & u \\ u & & \\ & u & 
                                           \end{pmatrix} .\]
Hence $\w_A$ cannot be a free $A$-module of rank one and $A$ is not Gorenstein. 
However, if we define $b'$ by $b'_\s=b'_\t=b_\s$ then the corresponding order has the same ramification data and is Gorenstein.

\end{subsection}

\begin{subsection}{Conclusion}
We conclude the following theorems from
the above constructions. 
\begin{theorem}   \label{gor} 
A canonical order is Gorenstein. 
\end{theorem}
We also observe the following simple version of the McKay correspondence.
\begin{theorem}
For a canonical order singularity, the number of permissible modules equals $n+1$, 
where $n$ is the number of exceptional curves of the minimal resolution.
\end{theorem} 
\end{subsection}

Finally, given the $G'$-modules computed, it is easy to compute the AR quivers of all the canonical orders using Prop.~\ref{arquivers}. 
They are listed in the appendix. 
We only remark the following interesting observations. 
First, a permissible module is either an irreducible $G'$-module or the direct sum of two such,
and the latter occurs precisely when the two modules are related by the Auslander-Reiten translation. 
Second, the AR quivers for the maximal canonical orders are precisely the quotient diagrams in \cite[Table~5.28]{A86}. 
This gives an intrinsic description of these diagrams. 

\end{section}

\begin{section}{Appendix}
\subsection{Minimal resolutions}

We list the intersection graphs of the exceptional and ramification curves for the minimal resolutions of canonical orders. 
These were determined in Thm.~\ref{canramdata}. 
The vertices of types $\bullet$, $\circ$, and $\odot$ denote exceptional curves, ramification curves, and exceptional ramification curves
respectively.
Two vertices are joined by a single edge (resp. a double edge) if the corresponding curves intersect transversely (resp. tangentially)
in a single point. 
Each exceptional curve is labelled above by the negative of its self-intersection.
We also label ramification curves below by the ramification index $e$ when $e \neq 2$.  
We do not include the well known Dynkin diagrams for types $A_n$, $D_n$, and $E_n$.

\begin{figure}[hp] \centering \caption{Minimal resolutions of canonical orders}
\begin{eqnarray*}
\begin{array}{lc}
A_{1,2,\xi} &
\xymatrix{ 
\overset{\phantom{1}}{\underset{2e}{\circ}}
\ar@{-} [r] & 
\overset{1}{{\underset{e}{\odot}}}
\ar@{-} [r] & 
\overset{\phantom{1}}{\underset{2e}{\circ}}
}
\\
\\
BL_n &
\xymatrix{  
\overset{2}{\underset{\phantom{1}}{\bullet}}
\ar @{-}[r] & 
\overset{2}{\underset{\phantom{1}}{\bullet}} 
\ar @{.}[r] & 
\overset{2}{\underset{\phantom{1}}{\bullet}} 
\ar @{-}[r] & 
\overset{1}{\underset{\phantom{1}}{\bullet}} 
\ar @{=}[r] & 
\circ  
}
\\
\\
B_n &
\xymatrix{   
& & & & \circ \\
\overset{2}{\underset{\phantom{1}}{\bullet}}
\ar @{-} [r] & 
\overset{2}{\underset{\phantom{1}}{\bullet}}
\ar @{.} [r] & 
\overset{2}{\underset{\phantom{1}}{\bullet}}
\ar @{-} [r] & 
\overset{1}{\underset{\phantom{1}}{\bullet}} 
\ar @{-} [ur]  \ar @{-} [dr] \\
  & & & &  {\circ}
}
\\
\\ 
L_n &
\xymatrix{  
& & & & \circ \\
\overset{2}{\underset{\phantom{1}}{\bullet}}
\ar @{-} [r] & 
\overset{2}{\underset{\phantom{1}}{\bullet}}
\ar @{.} [r] & 
\overset{2}{\underset{\phantom{1}}{\bullet}}
\ar @{-} [r] & 
\overset{1}{\underset{\phantom{1}}{\bullet}} 
\ar @{-} [ur]  \ar @{-} [dr] \\
  & & &  & \circ
\ar @{-} [uu]  
}
\\
\\
DL_n &
\xymatrix{
\circ
\ar @{-}[r] & 
\overset{2}{{\odot}} 
\ar @{-}[r] & 
\overset{2}{{\odot}}
\ar @{.}[r] & 
\overset{2}{{\odot}}
\ar @{-}[r] & 
\overset{2}{{\odot}}
\ar @{-}[r]  \ar @{-} [d] &
\overset{1}{{\bullet}} \ar@{-}[dl] \\
& & &  & \circ &
} 
\\
\\
BD_n &
\xymatrix{  
{ \circ}
\ar @{-} [r] & 
\overset{2}{{\odot}}
\ar @{-} [r] & 
\overset{2}{{\odot}}
\ar @{.} [r] & 
\overset{2}{{\odot}}
\ar@{-} [r] &
\overset{2}{{\odot}}
\ar @{-} [r]&  
\overset{1}{{\bullet}}
\ar@{-} [r] &  
{ \circ } \\
& & & & { \circ } \ar@{-}[u]&  & 
}
\\
\\
A_{n,\xi} & 
\xymatrix{ 
\overset{\phantom{1}}{\underset{e}{\circ}}
\ar @{-} [r] & 
\overset{2}{\underset{e}{\odot}}
\ar @{-} [r] & 
\overset{2}{\underset{e}{\odot}}
\ar @{.} [r] & 
\overset{2}{\underset{e}{\odot}}
\ar@{-} [r] &
\overset{2}{\underset{e}{\odot}}
\ar @{-} [r] &  
\overset{\phantom{1}}{\underset{e}{\circ}}
}
\end{array}
\end{eqnarray*}
\end{figure}
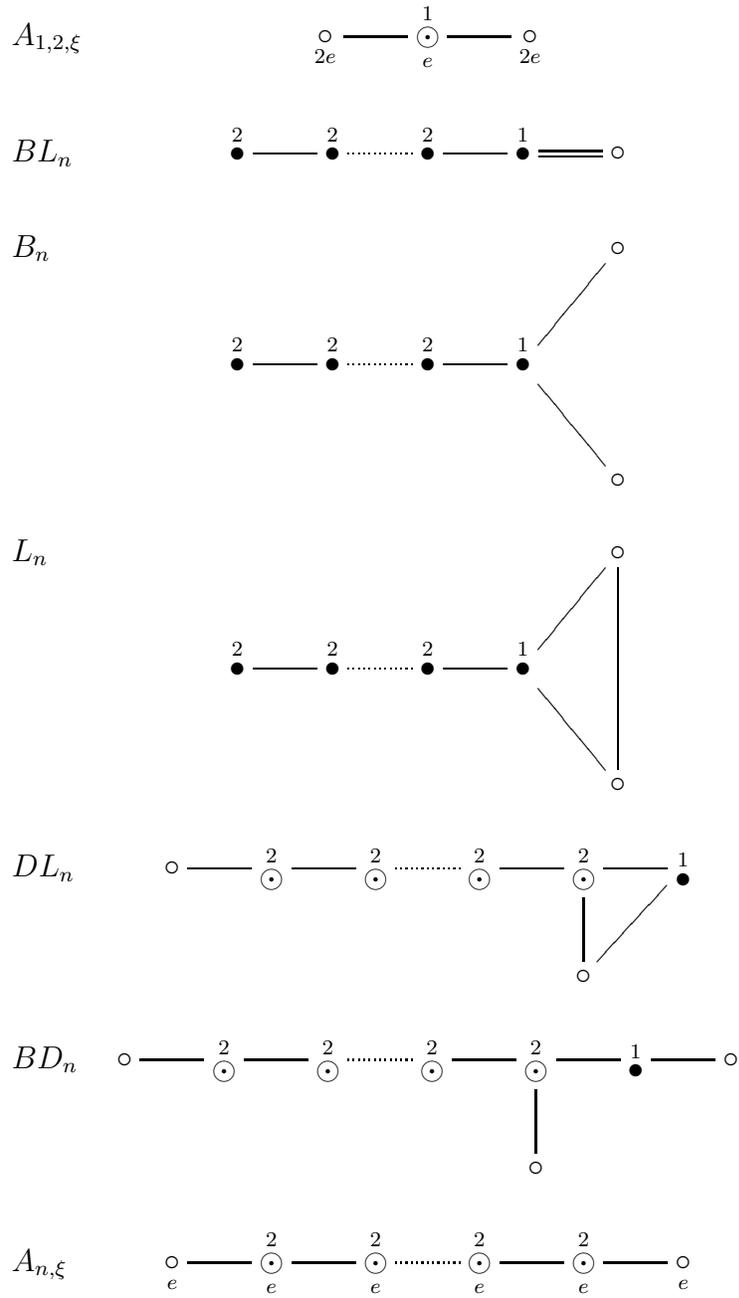

\clearpage

\subsection{Auslander--Reiten quivers}

We list the AR quivers for canonical orders $A$.
They are obtained from Prop.~\ref{arquivers} and the explicit computations in Sec.~\ref{exp}. 
The vertices correspond to isomorphism classes of indecomposable reflexive $A$-modules.
The solid arrows correspond to morphisms of $A$-modules.
The dashed arrows give the action of the AR translation $\tau$.
See \cite[Def.~5.2]{Yoshino} for more details.
We do not include the affine Dynkin diagrams for types $A_n$, $D_n$, and $E_n$.

\begin{figure}[hp] \centering \caption{Auslander--Reiten quivers of canonical orders I}
\begin{eqnarray*}
\begin{array}{lc}
A_{1,2,\xi} &
\xymatrix{
1 \ar@<.25ex> [r] \ar@<.25ex> [d] \ar@<.25ex>@{.>} [rd] & 2 \ar@<.25ex> [l] \ar@<.25ex> [d] \ar@<.25ex>@{.>} [ld] \\
3 \ar@<.25ex> [u] \ar@<.25ex> [r] \ar@<.25ex>@{.>} [ru] & 4 \ar@<.25ex> [l] \ar@<.25ex> [u] \ar@<.25ex>@{.>} [ul]
}
\\
\\
BL_n &
\xymatrix{
1  \ar@<.25ex>@{.>} [dd] \ar@<.25ex> [dr] & {} & {} &  &  {} \\
{}  & 3  \ar@{.>}@(ul,ur) \ar@<.25ex> [r] \ar@<.25ex> [ul] \ar@<.25ex> [ld] & 4  \ar@{.>}@(ul,ur) \ar@<.25ex> [l] \ar@{.} [r]
& n  \ar@{.>}@(ur,ul) \ar@<.25ex> [r] & n+1 \ar@(dl,dr) \ar@{.>}@(ur,ul) \ar@<.25ex> [l]  \\
2 \ar@<.25ex> @{.>} [uu] \ar@<.25ex> [ur] & {} & &  {} & {} 
}
\\
\\
B_n &
\xymatrix{
1  \ar@<.25ex>@{.>} [dd] \ar@<.25ex> [dr] & {} & {} &  &  {} & 
n+2 \ar@<.25ex> @{.>} [dd] \ar@<.25ex> [dl] \\
{}  & 3   \ar@{.>}@(ul,ur) \ar@<.25ex> [r] \ar@<.25ex> [ul] \ar@<.25ex> [ld] & 4  \ar@{.>}@(ul,ur) \ar@<.25ex> [l] \ar@{.} [r]
& n  \ar@{.>}@(ur,ul) \ar@<.25ex> [r] & n+1 \ar@{.>}@(ur,ul) \ar@<.25ex> [l] \ar@<.25ex> [ur] \ar@<.25ex> [dr] \\
2 \ar@<.25ex> @{.>} [uu] \ar@<.25ex> [ur] & {} & &  {} & {} & 
n+3 \ar@<.25ex>@{.>} [uu]  \ar@<.25ex> [ul] 
}
\\
\\
L_n &
\xymatrix{ 
\\ 
1  \ar@(dr,ld) \ar@{.>}@(ur,ul)
\ar@<.4ex>[r] & 
2  \ar@{.>}@(ur,ul)
\ar@<.4ex>[r]
\ar@<.4ex>[l] & 
 3   
\ar@{.}[r] \ar@{.>}@(ur,ul)
\ar@<.4ex>[l]&
 n  \ar@{.>}@(ur,ul)
\ar@<.4ex>[r] &  
n+1  \ar@{.>}@(ur,ul) \ar@(rd,dl) \ar@<.4ex>[l] 
}
\end{array}
\end{eqnarray*}
\end{figure}

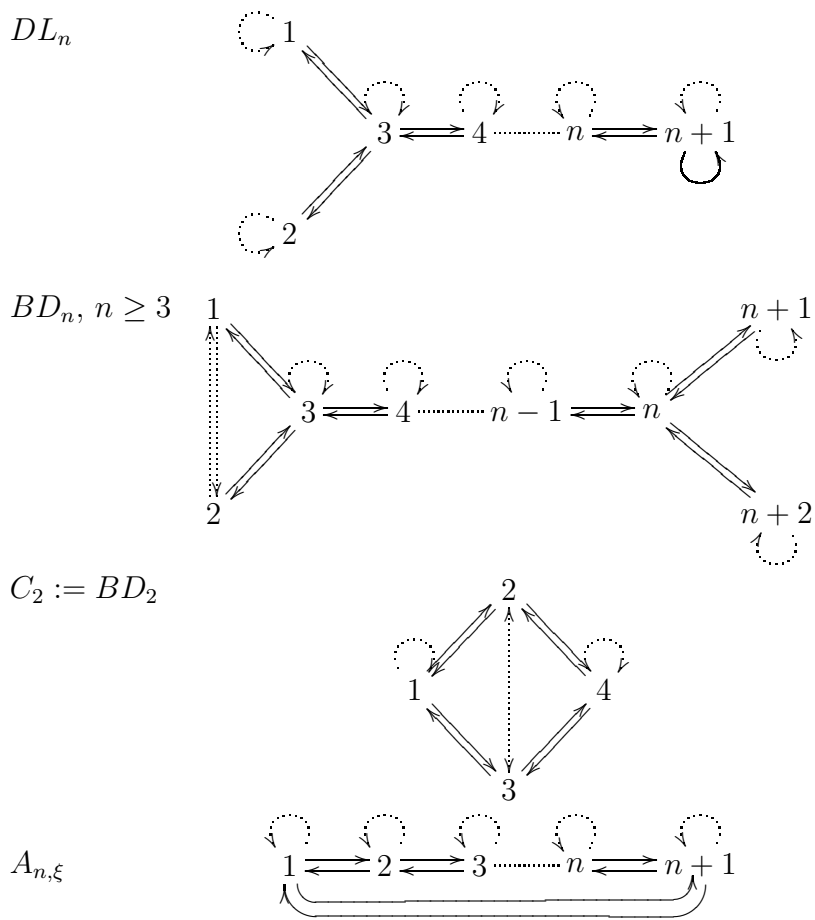
\begin{figure}[hp] \centering \caption{Auslander--Reiten quivers of canonical orders II}
\begin{eqnarray*}
\begin{array}{lc}
DL_n &
\xymatrix{
1  \ar@{.>}@(ul,dl) \ar@<.25ex> [dr] & {} & {} &  &  {} \\
{}  & 3  \ar@{.>}@(ul,ur) \ar@<.25ex> [r] \ar@<.25ex> [ul] \ar@<.25ex> [ld] & 4  \ar@{.>}@(ul,ur) \ar@<.25ex> [l] \ar@{.} [r]
& n  \ar@{.>}@(ur,ul) \ar@<.25ex> [r] & n+1 \ar@(dl,dr) \ar@{.>}@(ur,ul) \ar@<.25ex> [l]  \\
2  \ar@{.>}@(ul,dl) \ar@<.25ex> [ur] & {} & &  {} & {} 
}
\\
\\
BD_n, \, n \ge 3 &
\xymatrix{
1  \ar@<.25ex>@{.>} [dd] \ar@<.25ex> [dr] & {} & {} &  &  {} & 
n+1 \ar@(dl,dr)@{.>}\ar@<.25ex> [dl] \\
{}  & 3  \ar@{.>}@(ul,ur) \ar@<.25ex> [r] \ar@<.25ex> [ul] \ar@<.25ex> [ld] & 4  \ar@{.>}@(ul,ur) \ar@<.25ex> [l] \ar@{.} [r]
& n-1  \ar@{.>}@(ur,ul) \ar@<.25ex> [r] & n \ar@{.>}@(ur,ul) \ar@<.25ex> [l] \ar@<.25ex> [ur] \ar@<.25ex> [dr] \\
2 \ar@<.25ex> @{.>} [uu] \ar@<.25ex> [ur] & {} & &  {} & {} & 
n+2 \ar@(dr,dl)@{.>} \ar@<.25ex> [ul] 
}
\\
\\
C_2:=BD_2 &
\xymatrix{
{} & 2 \ar@{.>} [dd] \ar@<.25ex> [ld] \ar@<.25ex> [dr] & {} \\
 1  \ar@{.>}@(ul,ur) \ar@<.25ex> [ur] \ar@<.25ex> [dr] & {} & 4  
\ar@{.>}@(ul,ur) \ar@<.25ex> [lu] \ar@<.25ex> [ld]  \\
{} & 3  \ar@{.>} [uu]  \ar@<.25ex> [lu] \ar@<.25ex> [ur] & {} 
}
\\
\\
A_{n,\xi} &
\xymatrix{ 
1  \ar@<.4ex> `d[r] `[rrrr] \ar@{.>}@(ur,ul)
\ar@<.4ex>[r] & 
2  \ar@{.>}@(ur,ul)
\ar@<.4ex>[r]
\ar@<.4ex>[l] & 
 3   
\ar@{.}[r] \ar@{.>}@(ur,ul)
\ar@<.4ex>[l]&
 n  \ar@{.>}@(ur,ul)
\ar@<.4ex>[r] &  
\ar@<.4ex>[l]
n+1  \ar@{.>}@(ur,ul) \ar@<.4ex>@{-<} `d[l] `[llll]
}
\end{array}
\end{eqnarray*}
\end{figure}
\clearpage

\end{section}

\medskip
\noindent
Daniel Chan, School of Mathematics, University of New South Wales, Sydney, 2052, NSW, Australia.\\
e-mail: \texttt{danielch@unsw.edu.au} \\ 
\\
Paul Hacking,  Department of Mathematics, University of Washington, Box 354350, Seattle, WA~98195, USA.\\ 
e-mail: \texttt{hacking@math.washington.edu} \\
\\
Colin Ingalls, Department of Mathematics and Statistics, University of New Brunswick, Fredericton, NB E3B 5A3, Canada.\\
e-mail: \texttt{colin@math.unb.ca} \\

\end{document}